\documentclass[graybox]{svmult}

\usepackage{mathptmx}       
\usepackage{helvet}         
\usepackage{courier}        
\usepackage{type1cm}        
\usepackage{amssymb,amsmath,latexsym} 
\numberwithin{equation}{section}               
\usepackage{makeidx}         
\usepackage{graphicx}        
\graphicspath{{Figures/}}
\usepackage{multicol}        
\usepackage[bottom]{footmisc}

\makeindex             

\begin{document}

\title*{ Quasi-linear elliptic equations with data in $L^{1}$ on a compact Riemannian manifold}
\author{E. AZROUL, A. ABNOUNE and M.T.K. ABBASSI}
\institute{E. AZROUL \at  Laboratory of Mathematical Analysis and Applications,
	Faculty of sciences Dhar El Mahraz,
	Sidi Mohamed Ben Abdellah University,Fes, Morocco, \email{elhoussine.azroul@gmail.com}
\and A. ABNOUNE \at  Laboratory of Mathematical Analysis and Applications,
	Faculty of sciences Dhar El Mahraz,
	Sidi Mohamed Ben Abdellah University,Fes, Morocco, \email{abdelabn@gmail.com}
	\and M.T.K. ABBASSI \at  Laboratory of Mathematical ALgebra and geometry,
	Faculty of sciences Dhar El Mahraz,
	Sidi Mohamed Ben Abdellah University,Fes, Morocco, \email{mtk$_{-}$abbassi@yahoo.fr}}


\maketitle

\abstract{This work is dedicated to the study of quasi-linear elliptic problems with $L^1$ data, the simple model will be the next equation on $ (M,g) $ a compact Riemannian manifold. 
 $$-\Delta_{p} u=f$$
 Where $f\in L^{1}(M) $ .Our goal is to develop the functional framework and tools that are necessary to prove the existence and the uniqueness of the solution for the previous problem. Notice that our argument can be used to deal with a more general class of quasi-linear equations.\\}

\textbf{Introduction}\\

This article is dedicated to the study of quasi-linear elliptic equations with data in $ L ^ {1}(M) $, the major difficulty encountered when one is interested in such problems is that the classical theories of existence, either using variational methods or compactness methods, are not applicable. Hence the need to use new techniques to prove the existence and uniqueness of solutions for such problems.

Note that the importance of trying to solve problems with data in $ L ^ {1} $ is not limited to a purely theoretical framework, but also for applicable reasons, to be convinced it is recommended to the reader various references as for example: $ [1], [2] $ and $ [3] $ where different examples of equations having an application in physics are presented.\\

The first significant advance in this direction is due to Stampacchia in $ [3], $ where he considers second-order linear elliptic operators with non-regular data of the form $$ L (u) = f $$ where
$$L(u)=\sum_{i, j=1}^{N} a_{i j} \frac{\partial^{2} u}{\partial x_{i j}}+\sum_{i=1}^{N} \frac{\partial u}{\partial x_{i}}+c u$$
where $a_{i j}, b_{i}, c$ are functions with specific hypotheses .

In his famous works, Stampacchia uses the notion of duality to solve these classes of problems. Existence and uniqueness results have been proved in this direction thanks to the linear character and the "regularizing effect" of the operator.
Note that in the case where $L \equiv \Delta$ 
then the notion of duality coincides with the notion of solution in the sense of distributions. In the linear framework and with the notion of duality,we can even consider data measures.

The extension of Stampacchia's work to non-linear operators has been done by several mathematicians. The first works were realized by Boccardo, Mu-rat, Gallouet and their collaborators. The main difficulties for non-operators
linear consists of two points:
\begin{enumerate}
\item The sense in which the solution is defined (the meaning of the good solution and the method of its construction).
\item The uniqueness of the "good" solution.
\end{enumerate}

Note that the second question is legitimate given Serrin's counter example for the non-uniqueness of the solution, see $[4]$.

To go beyond the first difficulty we proceed by approximation by returning
In the variational framework, the main step is to demonstrate properties of the solutions for approximate problems that remain conserved by passing to the limit.
This passage is feasible by imposing natural conditions on the space of the test functions.

concerning the second difficulty,we demonstrate partial results, especially
 for the $ \Delta_ {p} $ operator we are able to demonstrate the uniqueness of the solution. It will be noted that the uniqueness of the solution is usually true.\\
 
 We organize this work in two sections.In the first section we briefly recall the functional spaces of Sobolev and Marcinkiewicz, on a compact Riemannian manifold, which will be very useful in this paper.In Subsection 1.2.3 we define the notion of the weak solution. Using variational techniques we prove the existence and the uniqueness of the energy solution for the problem :
$$-\Delta_{p} u=f, \quad u \in W_{0}^{1, p}(M)$$ where $ (M,g) $  a  Riemannian manifold.This result is a natural extension of the Lax-Milgram Theorem to the non-linear case in a Euclidean space. In Section 1.3 we present the proof of the Picone inequality in its general version on a compact Riemannian manifold and as a consequence we obtain a comparison principle for quasi-linear problems with a "concave" term compared to Laplacian.This result generalizes that of Brezis-Kamin in $ [5] $ for the Laplacian, see $ [6] $.The second section is dedicated to define the notion of entropy on a compact Riemannian manifold, in which we will study our problem. We begin by defining the functional framework that will be given using the truncation function, ie we analyze the functional properties of $ T_{k}(u)$ instead of $ u.$ Note that $ T_ {k}: \mathbb {R} \rightarrow \mathbb {R} $ defined by :

$$T_{k}(s)=\left\{\begin{array}{ccc}{s} & {\text { if }} & {|s| \leq k} \\ {k \operatorname{sign}(s)} & {\text { if }} & {|s|>k}\end{array}\right.$$

After giving the definition of solution in the sense of entropy, we prove the existence and uniqueness of the solution in this context, some properties of the entropy solution in Marcinkiewicz's spaces on Riemannian manifolds will be deduced. At the end of the section
some generalizations for non-homogeneous quasi-linears operators and with second members that may depend on $ u $ will be presented, see $ [1] $.

\section{Preliminaries}

\subsection{some definitions}

\begin{definition}[Equi-integrable functions in $L^1$]

Let X be a set of $ \mathbb {R} ^ {N} $.
We say that a sequence $ \left \{f_ {n} \right \} $ of functions of $ L ^ {1} (X) $ is equi-integrable if, for all $ \varepsilon> 0 $, there exists $ \delta> 0 $ such as meas $ (E) <\delta $ with $ E \subset X $ will result for all $ n $,
$$ 
\int_{E}\left|f_{n}(x)\right| d x \leq \varepsilon.$$
\end{definition}

We use
often the next result of compactness in $L^{1}$.
\begin{lemma}[Lemma of Vitali : compactness in $L^1$]

Let $X$ a finite measure set for the Lebesgue measure of $\mathbb{R}^{N} .$ Let $\left\{f_{n}\right\}$ a sequence of functions of $ L ^ {1} (X) $ which converges everywhere to $ f, $ and which is equi-integrable. Then $ f \in L ^ {1} (X) $ and $ \left \{f_ {n} \right \} $ converges strongly to $ f $ in $ L ^ {1} (X) $.
\end{lemma}

\subsection{Functional spaces}
\subsubsection{Sobolev spaces $ W ^ {k, p} (M) $}
see$ [7] $ and $ [8] $

Let $(M,g)$ a Riemannian manifold, for an integer $k$ and $u\in C{^{\infty}}(M)$ , $\nabla ^{k}u$ represents the $k-th$ of the covariant derivative of $u$ (with the Convention
$\nabla ^{0} {u} = u$ ) . and the norm of $k-th$ of covariant derivative on a local map is given by the  formula :

$$|\nabla ^{k} {u}|= g^{{i_{1}} {j_{1}}}.......g^{{i_{k}} {j_{k}}} ({\nabla ^{k} {u})}_{{{i_{1}....i_{k}}}} ({\nabla ^{k} {u})}_{{{j_{1}....j_{k}}}}$$ where the Einstein summation convention is adopted.

We also recall the notion of Riemannian measureon manifolds , let $\{{U_{i}} ,{\Phi_{i}}\}$ be any atlas of M . There exists a partition of unity $\{{U_{i}} ,{\Phi_{i}},{\eta_{i}}\}$ subordinate to $\{{U_{i}} ,{\Phi_{i}}\}$  Give a continuous function $f : M \rightarrow \mathbb{R}$ we define the integrale as follows 
$$ 
\int_{M}fd\sigma_{g}=\sum_{i} \int_{\Phi(U_{i})} ({\eta_{i}}\sqrt{det g} f\circ \Phi_{i}) dx.$$
where $dx$ is the Lebesgue measure on $ \mathbb{R}^{n}.$\\

Let be $p\geq  1$ a real, and $k$ a positive integer.

$$L^{p}{(M)}=\{ u:M\to \mathbb {R} \quad measurable \ / \int_{M} {|u|^{p} d\sigma_{g}<\infty\}}$$

\qquad 	$C_{k} ^{p} (M)$ functions space $u \in C^{\infty }$ such as
${|\nabla^{j} {u}| \in L^{p}{(M)} }$ for
$j=0,....,k$

$$C_{k} ^{p} (M) = \{ u \in C^{\infty} \ /  \forall j=0,....,k \quad  \int_{M} {|\nabla^{j} {u}|^{p} d\sigma_{g}<\infty\}}$$
\begin{definition}
 The Sobolev space $W^{k,p}(M)$ is the complete space $C^{p}_{k}(M)$ for the norm
$$\|u\|_{W^{k,p}(M)}=\sum_{j=0}^{k}\|\nabla^{j}u\|_{L^{p}(M)}$$

$$\|u\|_{W^{1,p}(M)} =\|\nabla u\|_{p} + \|u\|_{p} $$
\end{definition}
\begin{definition}
	
We must recall the notion of the geodesic distance for every curve : $$\varUpsilon : \ [ a,b  \ ] \to M $$
We define the length of $\varUpsilon$ by :
$$ l(\varUpsilon)  = \int _{a} ^{b} \sqrt { g(\varUpsilon(t)) ( \frac {d \varUpsilon}  {dt},\frac {d \varUpsilon}  {dt} )} dt$$
\end{definition}
\begin{remark}

For
$x,y\in M$ defining a distance $d_{g}$ by :

$$d_{g} (x,y) = inf \{l(\varUpsilon) : \varUpsilon : [0,1]  \to M \quad . \quad \varUpsilon(0)=x \quad, \quad \varUpsilon(1)=y        \}$$

By the theorem of Hopf-Rinow, we obtain that if $ M $ a Riemannian manifold then compact for all $ x, y $ in $ M $ can be joined by a minimizing curve $\varUpsilon $  i.e $l(\varUpsilon) = d_{g} (x,y) $
\end{remark}

\begin{proposition}
	If $p=2$,
	$W^{k,2}(M)$ is a Hilbert space 
space for the scalar product
	$$\left(u,v\right)_{H^{k}}=\sum_{j=0}^{k}\left(\nabla^{j}u,\nabla^{j}v\right)_{L^{2}}$$.
\end{proposition}

\begin{proposition}
	If $ p> 1 $ then $ W ^ {k,p} (M) $ is reflexive.
\end{proposition}
\begin{proposition}
	Any reflex normalized space is a Banach space.
	Then if $ p> 1 $ then $ W ^ {k,p} (M) $ is Banach.
\end{proposition}

\begin{definition}
The Sobolev space $ W_ {0} ^ {k, p} (M) $ is the closure of $\mathcal{D}(M)$ in $W^{k,p}(M)$ .
\end{definition}

\begin{theorem}
	If $(M,g)$ is complete, then for all $p\geq1$ $W^{1,p}_{0}(M)=W^{1,p}(M)$.\\
\end{theorem}

\textbf{Embeddings of Sobolev:}See $[7]$.\\

\begin{lemma}
	Let $(M,g)$ a complete Riemannian manifold of dimension $ n $. Suppose that inclusion $W^{1,1}(M)\subset L^{\frac{n}{(n-1)}}(M)$ is valid. Then, for a whole real number $1\leq q<p$ and an integer $0\leq m<k$ which verify $\frac{1}{p}=\frac{1}{q}-\frac{(k-m)}{n}$, $W^{k,q}(M)\subset W^{m,p}(M)$.
\end{lemma}

\begin{remark}

Note that the proof of the Lemma \ref{lem1} shows that if $A\in\mathbb{R}$ is such that $\forall~~u\in W^{1,1}(M)$,
$$\left(\int_{M}|u|^{n/(n-1)}d\sigma_{g}\right)^{(n-1)/n}\leq A\int_{M}\left(|\nabla u|+|u|\right)d\sigma_{g})$$
So, for all $1\leq q<n$ and all
$u\in W^{1,q}(M)$,
$$\left(\int_{M}|u|^{p}d\sigma_{g}\right)^{1/p}\leq \frac{A p(n-1)}{n}\left\{\left(\int_{M}|\nabla u|^{q}d\sigma_{g}\right)^{1/q}+\left(\int_{M}|u|^{q}d\sigma_{g}\right)^{1/q}\right\}$$
Where
$1/p=1/q-1/n.$

\end{remark}

\begin{theorem}
	Let $(M,g)$ a compact Riemannian manifold of dimension $ n $. For a real number $1\leq q<p$ and an integer $0\leq m<k$ which verify $\frac{1}{p}=\frac{1}{q}-\frac{(k-m)}{n}$, $W^{k,q}(M)\subset W^{m,p}(M)$.
\end{theorem}
\begin{theorem}\textbf{(Rellich-Kondrakov's Theorem):}
	Let $(M,g)$ a compact Riemannian manifold of $n$ dimension , $j\geq0$ and $m\geq1$ two integers, $q\geq1$ and $p$ two real numbers that verify $1\leq p<nq/(n-mq)$,  the inclusion $W^{j+m,q}(M)\subset W^{j,p}(M)$ is compact
\end{theorem}
\begin{corollary}
	Let $(M,g)$ a compact Riemannian manifold of $ n $ dimension.  For everything $1\leq q<n$ and  $p\geq1$ such as
	$\frac{1}{p}>\frac{1}{q}-\frac{1}{n}$, the inclusion $W^{1,q}(M)\subset L^{p}(M)$ is compact.
\end{corollary}
\begin{lemma}\label{lem4}\textbf{(Inequality of Poincare):}
	Let $D$ a regular domain is bounded in a Riemannian manifold $M$ and $1\leq p<\infty$. Then there is a constant $A$ such as:
	$$\left(\int_{D}|u-u_{D}|^{p}d\sigma_{g}\right)^{\frac{1}{p}}\leq A\left(\int_{D}|\nabla u|^{p}d\sigma_{g}\right)^{\frac{1}{p}},$$
	for everything
	$u\in W^{1,p}_{loc}(M)$, where  $u_{D}=\frac{1}{vol(D)}\int_{D}ud\sigma_{g}$ is the mean value of $ u $ on $ D $
\end{lemma}

By combining this lemma with the Holder inequality, we obtain:
\begin{corollary}\label{cor1}
	There exists a constant $ c = c_{D} $ such that
	\begin{align*}\label{(3)}
		\int_{D}|u-u_{D}|d\sigma_{g}\leq c_{D}\left(\int_{M}|\nabla u|^{p}d\sigma_{g}\right)^{\frac{1}{p}}~~~~~~\forall~u\in W^{1,p}_{loc}(M)
		\end{align*}
\end{corollary}


\subsubsection{Marcinkiewicz's spaces.}
\begin{definition}$ [7] $
Let $ f: M \rightarrow \mathbb {R} $ be a measurable function, its distribution function
$$ \phi_{f}(k) = \operatorname {meas}\bigg \{x \in M :|f(x)|>k\bigg\} \quad k>0
, $$
\end{definition}

\begin{definition}

 Let $ 0 <q <\infty $ and $ (M, g) $  Riemannian manifold, the space Marcinkiewicz $ \mathcal {M} ^ {q} (M) $ is the set of functions measurable $ f: M \rightarrow \mathbb {R} $ such as
$$ \phi_ {f} (k) \leq C k ^ {- q}, ~~ C <\infty, $$
\end{definition}

Marcinkiewicz's space $\mathcal{M}^{q}(M)$  is defined the norm
$$\|f\|_{\mathcal{M}^{q}(M)}=\inf \Bigg\{C : \phi_{f}(k) \leq C k^{-q}, \quad \text{for all}\quad k>0\Bigg\}.$$
is a Banach space.

Note that if $f \in L^{q}(M),$ we have
$$\int_{\{|f|>k\}} d\sigma_{g} \leq \int_{M}\left|\frac{f}{k}\right|^{q} d\sigma_{g} \leq k^{-q} \int_{M}|f|^{q} d\sigma_{g},$$

so
$$ 
\phi_{f}(k) \leq k^{-q}\|f\|_{q}^{q}
 $$

and as a conclusion will have $L^{q}(M) \subset \mathcal{M}^{q}(M) .$

For analyze the properties of the spaces $ \mathcal {M} ^ {q} (M), $ needs some  the  next lemma 

\begin{lemma}
$If f \in L^{q}(M)$ so $$ 
\int_{M}|f(x)|^{q} d x=q \int_{0}^{+\infty} t^{q-1} \phi_{f}(t) d t.
 $$
\end{lemma}
\textbf{Proof.}

We start with the case where $q = 1$. Let
$$ 
H(t)=\left\{\begin{array}{lll}{1} & {\text { if }} & {t>0,} \\ {0} & {\text { if }} & {t<0,}\end{array}\right.
 $$
so
$$ 
H(|f(x)|-k)=\left\{\begin{array}{ll}{1} & {\text { if ~~}|f(x)|>k,} \\ {0} & {\text { if ~~}|f(x)|<k,}\end{array}\right.
 $$
so we have
$$ 
\begin{aligned} \int_{0}^{+\infty} \phi_{f}(k)  dk_{g} &=\int_{0}^{+\infty}\left[\int_{M} H(|f(x)|-k) d\sigma_{g}\right] dk_{g} \\ &=\int_{M}\left[\int_{0}^{+\infty} H(|f(x)|-k)  dk_{g}\right] d\sigma_{g}, \quad \text { ( Fubini) } \\
&= \int_{M}\left[\int_{\{|f(x)|>k\}} 1 dk_{g}\right] d\sigma_{g}=\int_{M}\left[\int_{0}^{|f(x)|} 1  dk_{g}\right] d\sigma_{g},\\
&=\int_{M}|f(x)| d\sigma_{g},
\end{aligned}
 $$
so $\displaystyle\int_{0}^{+\infty} \phi_{f}(k) dk_{g}=\displaystyle\int_{M}|f(x)| d\sigma_{g}$
and the result is demonstrated.\\

We now consider the general case $q>1 .$ Ask $g(x)=|f(x)|^{q}$ so, $g \in L^{1}(M)$
and
$$\phi_{g}(k)=\operatorname{meas}\bigg\{|g|>k\bigg\}=\operatorname{mes}\bigg\{|f|^{q}>k\bigg\}=\operatorname{meas}\bigg\{|f|>k^{\frac{1}{q}}\bigg\},$$
i.e. $\quad \phi_{g}(k)=\phi_{f}\left(k^{\frac{1}{q}}\right),$

so
$$ 
\int_{M}|g(x)| d\sigma_{g}=\int_{0}^{+\infty} \phi_{f}\left(k^{\frac{1}{q}}\right) dk_{g}
 ,$$

Ask $t=k^{\frac{1}{q}}$ so $k =t^{q} $ and $d k=q t^{q-1}, $ so that  $$\int_{M}|f(x)|^{q} d\sigma_{g}=q \int_{0}^{+\infty} t^{q-1} \phi_{f}(t) d t . \blacksquare$$

\begin{corollary}
If $q \in ] 1, \infty[,$ so
$$ 
L^{q}(M) \subset \mathcal{M}^{q}(M) \subset L^{q-\varepsilon}(M) \quad \forall \varepsilon > 0,
 $$
and for all $q, \hat{q} \in[1, \infty[$ we have 
$$\mathcal{M}^{q}(M) \subset \mathcal{M}^{\hat{q}}(M) \quad \text{if} \quad q \geq \hat{q}.$$
\end{corollary}

\textbf{Proof.}

Assuming that $f \in \mathcal{M}^{q}(M)$ and $\varepsilon>0,$ we have

$$ 
\begin{aligned} \int_{M}|f(x)|^{q-\varepsilon} d\sigma_{g} &=\int_{\{|f| \leq 1\}}|f(x)|^{q-\varepsilon} d\sigma_{g}+\int_{\{|f|>1\}}|f(x)|^{q-\varepsilon}d\sigma_{g}, \\ & \leq c_{1}+\int_{\{|f| \leq 1\}}|f(x)|^{q-\varepsilon} d\sigma_{g}, \\ & \leq c_{1}+\int_{M}|f(x)|^{q-\varepsilon} 1_{\{|f|>1\}} d\sigma_{g}, \\ & \leq c_{1}+\int_{M}|g(x)|^{q-\varepsilon} d\sigma_{g}, \end{aligned}
 $$
 
where $g(x)=|f(x)| 1_{\{|f|>1\}},$ so
 $$ 
\phi_{g}(t)=\operatorname{meas}\{|g|>t\} \leq \operatorname{meas}\left\{|f(x)| 1_{\{|f|>1\}}>t\right\}
 $$

implies that

$$ 
\begin{aligned} \int_{M}|f(x)|^{q-\varepsilon} d\sigma_{g} & \leq c_{1}+\int_{M}|g(x)|^{q-\varepsilon} d\sigma_{g} ,\\ & \leq c_{1}+(q-\varepsilon) \int_{0}^{+\infty} t^{q-\varepsilon-1} \phi_{g}(t) d t ,\\ & \leq c_{1}+(q-\varepsilon) \int_{0}^{1} t^{q-\varepsilon-1} \phi_{g}(t) d t+(q-\varepsilon) \int_{1}^{+\infty} t^{q-\varepsilon-1} \phi_{g}(t) d t ,\\ & \leq c_{1}+c_{2}(q-\varepsilon) \int_{0}^{1} t^{q-\varepsilon-1} d t+c_{3}(q-\varepsilon) \int_{1}^{+\infty} t^{q-\varepsilon-1} t^{-q} d t ,\\
&  \leq c_{1}+c_{2}(q-\varepsilon)+c_{3}(q-\varepsilon) \int_{1}^{+\infty} t^{-\varepsilon-1} d t,\\
& \leq c_{4}+(q-\varepsilon)\left[-\frac{t^{-\varepsilon}}{\varepsilon}\right]_{1}^{+\infty},\\
&\leq C<\infty.
\end{aligned}
 $$

so $\quad f \in L^{q-\varepsilon}(M),$ and then it results that $\mathcal{M}^{q}(M) \subset L^{q-\varepsilon}(M) .$
As $L^{q}(M) \subset \mathcal{M}^{q}(M) \subset L^{q-\varepsilon}(M) \subset \mathcal{M}^{q-\varepsilon}(M) \quad$ for all $q \in ] 1, \infty[$ and for all
$\varepsilon>0,$ so we deduct for all $q, \hat{q} \in ] 1, \infty[$ wa have
$$\mathcal{M}^{q}(M) \subset \mathcal{M}^{\hat{q}}(M) \quad si \quad q \geq \hat{q}.~~\blacksquare$$

\subsubsection{Elliptic problems and the concept of the weak solution.} 

Let $ p> 1 $ and $ (M,g) $  a compact Riemannian manifold  , for $ u \in W_ {0} ^ {1, p} (M), $ we can consider the
continuous linear form $ - \Delta_ {p} u $ over $ W_ {0} ^ {1, p} (M) $ defined by
$$ \left \langle- \Delta_ {p} u, v \right \rangle \equiv \int _ {M} | \nabla u | ^ {p-2} \nabla u \nabla v d\sigma_{g}. $$
 
It's clear that $-\Delta_{p} u \in\left(W_{0}^{1, p}(M)\right)^{\prime}=W_{0}^{-1, p^{\prime}}(M)$ and $\left\|-\Delta_{p} u\right\|_{W_{0}^{-1, p^{\prime}}(M)}=\|u\|_{W_{0}^{1, p}(M)}$.

As a consequence, we have the next definition
 
\begin{definition}
let $f \in W_{0}^{-1, p^{\prime}}(m),$ we say that u is a weak solution of the problem
$$\left\{\begin{array}{ccc}{-\Delta_{p} u=f} & {\text {in }} & {M,} \\ {u(x)=0} & {\text { on }} & {\partial M,}\end{array}\right.$$

in $W_{0}^{1, p}(M)$ if and only if
$$\int_{M}|\nabla u|^{p-2} \nabla u \nabla \varphi d\sigma_{g}=\left\langle-\Delta_{p} u, \varphi\right\rangle \quad \forall \varphi \in W_{0}^{1, p}(M).$$
\end{definition}
 
The next inequalities will be systematically used in this work.

\begin{lemma}
Let $\xi_{1}, \xi_{2} \in M,$ we have
\begin{enumerate}
\item[1)] If $p\leq 2,$
$$\left|\xi_{1}+\xi_{2}\right|^{p}-\left|\xi_{1}\right|^{p}-p\left|\xi_{1}\right|^{p-2}\left\langle\xi_{1}, \xi_{2}\right\rangle \leq C(p)\left|\xi_{2}\right|^{p},$$

$$\left|\xi_{2}\right|^{p}-\left|\xi_{1}\right|^{p}-p\left|\xi_{1}\right|^{p-2}\left\langle\xi_{1}, \xi_{2}-\xi_{1}\right\rangle \geq C(p) \frac{\left|\xi_{2}-\xi_{1}\right|^{2}}{\left(\left|\xi_{2}\right|+\left|\xi_{1}\right|\right)^{2-p}}.$$

\item[2)] If $p > 2,$
$$\left|\xi_{1}+\xi_{2}\right|^{p}-\left|\xi_{1}\right|^{p}-p\left|\xi_{1}\right|^{p-2}\left\langle\xi_{1}, \xi_{2}\right\rangle \leq \frac{p(p-1)}{2}\left(\left|\xi_{1}\right|+\left|\xi_{2}\right|\right)^{p-2}\left|\xi_{2}\right|^{2},$$
$$\left|\xi_{2}\right|^{p}-\left|\xi_{1}\right|^{p}-p\left|\xi_{1}\right|^{p-2}\left\langle\xi_{1}, \xi_{2}-\xi_{1}\right\rangle \geq \frac{C(p)}{2^{p}-1}\left|\xi_{2}-\xi_{1}\right|^{p}.$$
\end{enumerate}
\end{lemma}

For the demonstration, see $[9]$ and $[10]$.\\

To demonstrate the existence of a weak solution for the previous problem, one often uses variational techniques and arguments of minimization of the convex functional ones. More precisely, we have the next result.
 
\begin{theorem}$[11]$
Let $ V $ be a reflexive Banach space, $ K \subset V $ is a closed non-empty convex, and $ J: K \rightarrow \mathbb {R} \cup \{+ \infty \} $ a semi-continue inferiorly coercive function 
 weakly on $ K. $

So $\inf _{u \in K} J(u)<\infty \quad$and $\quad \exists u_{0} \in K, \quad J\left(u_{0}\right)=\min _{u \in K} J(u).$

Moreover, if J is strictly convex, $ u_ {0} $ is unique. If $ J $ is differentiable in the sense of Gateaux and $ K $ then open $ J ^ {\prime} \left (u_ {0} \right) = 0. $
\end{theorem}

\begin{theorem}
Let$ (M,g) $  a  Riemannian manifold and $p \in ] 1, \infty[ .$ We
suppose that $f \in L^{q}(M)$ with $q \geq \overline{q}=\dfrac{N p}{N(p-1)+p},$ so there exists a unique solution  $u \in W_{0}^{1, p}(M)$ of the problem
$$ 
\left\{\begin{aligned}
-\Delta_{p} u=f & \quad\text { in }\quad&M , \\
 u=0 & \quad\text {if } & \partial M. \end{aligned}\right.
 $$
\end{theorem}

\subsection{Inequality of Picone for the $p-$ Laplacian and application.}

We begin by formulating the inequality of Picone punctual for the case of $p-$ Laplacian.
\begin{theorem}
Let $v>0, u \geq 0$ two positive class $C^{1}$ functions,  we pose
$$ 
L(u, v)=|\nabla u|^{p}+(p-1) \frac{u^{p}}{v^{p}}|\nabla v|^{p}-p \frac{u^{p-1}}{v^{p-1}}|\nabla v|^{p-2} \nabla v \nabla u.$$
$$R(u, v)=|\nabla u|^{p}-\nabla\left(\frac{u^{p}}{v^{p-1}}\right)|\nabla v|^{p-2} \nabla v.$$
so $L(u, v)=R(u, v), L(u, v) \geq 0$ and $L(u, v)=0,$ almost everywhere. in $M$   of $u=k v$ in each Connected component of $M .$
\end{theorem}
The proof of Theorem 6 is simple, it is based on the development of
term $ \nabla \left (\dfrac {u ^{p}}{v^{p-1}} \right) | \nabla v |^{p-2} \nabla v.$
To apply the Picone inequality to nonlinear elliptic equations
we need to prove an extension of Theorem $6$ in $W_{0}^{1,p} (M),$ more precisely we have the next lemma

\begin{lemma}
Let $v \in W^{1, p}(M)$ such as $v \geq \delta>0$ in $M .$ so for all $u \in C_{0}^{\infty}(M),$
$u \geq 0$
$$\int_{M}|\nabla u|^{p} \geq \int_{M}\left(\frac{|u|^{p}}{v^{p-1}}\right)\left(-\Delta_{p} v\right).$$
\end{lemma}

\textbf{Proof.}

As $v \in W^{1, p}(M) $ and $ v \geq \delta> 0 $ in $ M, $ then it exists a sequence $ \left \{v_ {n} \right \} $ regular functions such as
$$ 
\left\{\begin{array}{c}{v_{n} \rightarrow v \operatorname \quad{in}\quad W^{1, p}(M), v_{n} \in C^{1}(M)} \\ {v_{n} \rightarrow v\quad a.e ,\quad \text { and } v_{n}>\dfrac{\delta}{2}  \text \quad{in }\quad M .}\end{array}\right.
 $$

As a consequence of the continuity of the operator $-\Delta_{p}$ 
$\bigg(\operatorname{of}\quad W^{1, p}(M)$ in $W^{-1, p^{\prime}}(M),$ 
$p^{\prime}=\dfrac{p}{p-1} \bigg)$ we get that $-\Delta_{p} v_{n} \rightarrow-\Delta_{p} v$ in $W^{-1,~ p^{\prime}}(M), p^{\prime}=\dfrac{p}{p-1} \cdot(\operatorname{see}[24]) .$  En
using the identity of Picone at $ v_{n}, $ it results
$$ 
|\nabla u|^{p} \geq \nabla\left(\frac{u^{p}}{v_{n}^{p-1}}\right)\left|\nabla v_{n}\right|^{p-2} \nabla v_{n}.
 $$
as
$$ 
\begin{aligned} \int_{M}-\Delta_{p} v_{n} \frac{u^{p}}{v_{n}^{p-1}} &=\int_{M}\left|\nabla v_{n}\right|^{p-2}\left\langle\nabla v_{n}, \nabla\left(\frac{u^{p}}{v_{n}^{p-1}}\right)\right\rangle \\ &= p \int_{M} \frac{u^{p-1}}{v_{n}^{p-1}}\left|\nabla v_{n}\right|^{p-2}\left\langle\nabla v_{n}, \nabla u\right\rangle-(p-1) \int_{M} \frac{u^{p}}{v_{n}^{p}}\left|\nabla v_{n}\right|^{p} .\end{aligned}
 $$
 
Using the hypothesis on the Dominated convergence theorem we conclude
$$\int_{M}|\nabla u|^{p} \geq \int_{M}\left(\frac{-\Delta_{p} v}{v^{p-1}}\right) u^{p}, \quad u \in C_{0}^{\infty}(M), u \geq 0.$$ \hfill$\blacksquare$
 
In a more general context, we have the next result

\begin{theorem}
if $u \in W_{0}^{1, p}(M)$, $u \geq 0$,  $v \in W_{0}^{1, p}(M)$, $-\Delta_{p} v \geq 0$ is a measure of Radon bounded, $\left.v\right|_{\partial M}=0, v \gneq 0,$ so
$$\int_{M}|\nabla u|^{p} \geq \int_{M}\left(\frac{u^{p}}{v^{p-1}}\right)\left(-\Delta_{p} v\right).$$
\end{theorem}

\textbf{Proof. }  According to the principle of Maximum strong we have $ v>0 $ in $ M. $ (See $[13]) $.
We pose $ v_{m}(x) = v(x)+ \dfrac {1}{m}, m \in \mathbb {N}. $ So $ \Delta_ {p} v_{m} = \Delta_{p}v $ and $ \left \{v_ {m} \right \}$ converges in $W ^{1, p}(M) $ and a.e to $ v. $ Therefore, using Lemma $1.4,$ on gets the result for all $ \phi \in C_ {0} ^ {\infty} (M), \phi \geq 0.$ Now in the general case, by density we deduce
the existence of $u_{n} \rightarrow u$ in $W_{0}^{1, p}(M)$, $u_{n} \in C_{0}^{\infty}(M)$ et $u_{n} \geq 0$, so
$$ 
\int_{M}\left|\nabla u_{n}\right|^{p} \geq \int_{M}\left(\frac{-\Delta_{p} v_{n}}{v_{n}^{p-1}}\right) u_{n}^{p}=\int_{M}\left(\frac{-\Delta_{p} v}{v_{n}^{p-1}}\right) u_{n}^{p}.
 $$
By the hypothesis imposed on $u$ and according to the Lemma of Fatou we obtain the result.\hfill$\blacksquare$

\subsubsection{Comparison principle.}
As application of  lemma $6$, we demonstrate the next comparison result .
\begin{lemma}
Let $f$ be a continuous positive function such that $ \dfrac {f(u)} {u^{p-1}} \downarrow $ with $ 1<p. $
We
suppose that
  $ u,v \in W_{0}^{1, p}(M) \cap C ^{1} (M) $ are such that
$$ 
\left\{\begin{array}{ll}{-\Delta_{p} u \geq f(u),} & {u>0 \text { in } M} \\ {-\Delta_{p} v \leq f(v),} & {v>0 \text { in } M}\end{array}\right.
 $$
so $u \geq v$ in $M$
\end{lemma}
\textbf{Proof. } previous inequality implies that
$$ 
\frac{-\Delta_{p} u}{u^{p-1}}+\frac{\Delta_{p} v}{v^{p-1}} \geq \frac{f(u)}{u^{p-1}}-\frac{f(v)}{v^{p-1}}
 .$$
Multiply by $w=\left(v^{p}-u^{p}\right)^{+},$ we get that
$$ 
\begin{aligned} \int_{M}\left(\frac{-\Delta_{p} u}{u^{p-1}}+\frac{\Delta_{p} v}{v^{p-1}}\right)\left(v^{p}-u^{p}\right)^{+} & \geq \int_{M}\left(\frac{f(u)}{u^{p-1}}-\frac{f(v)}{v^{p-1}}\right)\left(v^{p}-u^{p}\right)^{+} \\ &=\int_{[v>-1]}\left(\frac{f(u)}{u^{p-1}}-\frac{f(v)}{v^{p-1}}\right)\left(v^{p}-u^{p}\right)^{+} \end{aligned}
 $$
 By the assumption on $f,$ we conclude that the term on the right in the previous equality is positive. On the other hand as $w=\left(v^{p}-u^{p}\right)^{+},$ so $\nabla w=p\left(v^{p-1} \nabla v-u^{p-1} \nabla u\right) \chi_{[v \geq u]}$,
so

$$ 
\begin{aligned} & \int_{M}\left(\frac{-\Delta_{p} u}{u^{p-1}}+\frac{\Delta_{p} v}{v^{p-1}}\right) w=\int_{M}|\nabla u|^{p-2}\left\langle\nabla u, \nabla\left(\frac{w}{u^{p-1}}\right)\right\rangle-\int_{M}|\nabla v|^{p-2}\left\langle\nabla v, \nabla\left(\frac{w}{v^{p-1}}\right)\right\rangle \\=& \int_{M}|\nabla u|^{p-2}\left\langle\nabla u, \frac{u^{p-1} \nabla w-(p-1) u^{p-2} w \nabla u}{u^{2(p-1)}}\right\rangle \\-& \int_{M}|\nabla v|^{p-2}\left\langle\nabla v, \frac{v^{p-1} \nabla w-(p-1) v^{p-2} w \nabla v}{v^{2(p-1)}}\right\rangle \\=& \int_{M \cap[v>u]}\left[\nabla\left.u\right|^{p-2}\langle\nabla u, \nabla v\rangle-(p-1) \frac{v^{p}}{u^{p}}|\nabla u|^{p}-|\nabla u|^{p}\right] \\ +& \int_{M \cap[v>u]}\left[p \frac{u^{p-1}}{v^{p-1}}|\nabla v|^{p-2}\langle\nabla v, \nabla u\rangle-(p-1) \frac{u^{p}}{v^{p}}|\nabla v|^{p}-|\nabla v|^{p}\right] \\=&
 \int_{M \cap[v>u]} K_{1}(x) d\sigma_{g}+\int_{M \cap[v>u]} K_{2}(x) d\sigma_{g}
\end{aligned}
 $$

and as $u>0$ and $v>0 $ in $ m $, using the Picone inequality, $K_{1} \leq 0$ and $K_{2} \leq 0 $. So $$\int_{m}\left(\frac{-\Delta_{p} u}{u^{p-1}}+\frac{\Delta_{p} v}{v^{p-1}}\right) w \leq 0 $$

and Consequently,
$$ 
\int_{m \cap[v \geq u]}\left(\frac{f(u)}{u^{p-1}}-\frac{f(v)}{v^{p-1}}\right)\left(v^{p}-u^{p}\right) \leq 0.
 $$

But on the set $[v>u]$, $
\dfrac{f(u)}{u^{p-1}}-\dfrac{f(v)}{v^{p-1}} \geq 0,
 $ so $|[v>u]|=0,$ and we deduce that $v \leq u .$
 
Easily demonstrates the extension
using Lemma $7$\hfill$\blacksquare$

\begin{lemma}[Comparison principle]
Let $u, v \in W_{0}^{1, p}(M) \cap C^{1}(M)$ such as
$$ 
\left\{\begin{array}{cc}{-\Delta_{p} u \geq h(x) f(u),} & {u>0 \text { in } M} \\ {-\Delta_{p} v \leq h(x) f(v),} & {v>0 \text { in } M}\end{array}\right.
 $$
\textbf{where $h$ is a positive function such that} $h \neq 0 .$ So $u \geq v$ in $M.$
\end{lemma}

\begin{remark}
The result of Lemma $ 8 $
is valid if $h(x)=|x|^{-p}$.
\end{remark}
As a direct application of Lemma $ 8, $ we obtain the next uniqueness result
\begin{theorem}
The problem
$$ 
\left\{\begin{array}{c}{-\Delta_{p} u=\lambda h(x) u^{q} \text { in } M, \quad 0<q<p-1} \\ {u>0 \text { in } M} \\ {\left.u\right|_{\partial M}=0}\end{array}\right.
 $$
where $h$ is in the conditions of the preceding theorem, admits a unique solution.
\end{theorem}

\begin{remark}
 In general, we have the same result of uniqueness if we replace $u ^{q}$ with a function of Carateodory $f(x,u)$ such that $ \dfrac {f(x,u)} {u^{p-1}}$ is
decreasing uniformly in $ x \in M.$ To demonstrate existence
we need to impose more conditions on $ f $.
\end{remark}

\section{Theory of existence and uniqueness of solutions for nonlinear elliptic problems with data in $L^{1}$}

\subsection{Introduction}
Consider the problem of form
\begin{equation}
 \left\{ {\begin{array}{*{20}{c}}
{ - {\Delta _p}u}& = &f&{{\rm{ in }}}&{M ,}\\
u& = &0&{{\rm{ on }}}&{\partial M .}
\end{array}} \right . \label{1}
\end{equation}

where $ 1<p<\infty,f $ is a measurable function such that $f \in L^{1}(M)$.

There are three difficulties associated with the study of the equation \eqref {1} .
\begin {enumerate}
\item [1-] Find the direction for which the previous equation is well defined.
\item [2-] The construction of a solution in the direction obtained.
\item [3-] Uniqueness of the solution found.
\end {enumerate}
Note that the most general meaning that can be used is the direction of distribution,
ie, u checks
$$\int_{M}|\nabla u|^{p-2} \nabla u \nabla \phi d\sigma_{g}=\int_{M} f \phi d\sigma_{g} \quad \forall \phi \in C_{0}^{\infty}(M)$$
except that the problem in this context is who
we do not have a construction argument
(the test function space being too "small"), and the second problem is the uniqueness of the solution ( the operator is nonlinear). Note that for the case $ p = 2, $ the distributional framework
  is a natural framework for studying equations with a second member in $ L^{1},$ because $ \Delta u=0 $ in the distributions sense implies that $u$ is harmonic in the classical sense.

To solve the nonlinear problem we need to introduce a new space $ \tau_{loc}^{1,1}(M) $
in which we can make sense of the gradient of $u,$ which in general is not
locally integrable. So the idea is to work with the truncations $ T_{k}(u)$ of the $u$ solution and expand the space of the test functions to bounded functions with a gradient in a suitable Lebesgue space.

The arguments we will introduce will be applicable to a class of equations general form.

$$ - {\Delta _p}u  = F(x,u)\quad in \quad{D^{ ' } (M)} \quad \quad (1,2)$$

Or $F$ is a carathéodory functions, continuous and decreasing in $ u $ for $ x $ fixed, and measurable in x for $ u $ fixed. moreover, $ F (x, 0) \in L^{1}(M) $ and $F(x,c) \in L_{loc}^{1}(M) $ if $c \neq 0 $ ,
and if
$$ 
G_{c}(x)=\sup _{|u| \leq c}|F(x, u)|,
 $$
  so $ G_{c} \in L_{l o c}^{1}(M) $ for all $c>0 $
 
\subsection{Functional Framework}
Before discussing the concept of the entropy solution, we will present the functional framework in which the solution is well defined.
We start with the introduction of the truncation operator. For a constant $k>0,$ we define the function $T_{k}: \mathbb{R} \rightarrow \mathbb{R}$ by

$$ 
T_{k}(s)=\left\{\begin{array}{ccc}{s} & {\text { if }} & {|s| \leq k,} \\ {k \operatorname{sign}(s)} & {\text { if }} & {|s|>k.}\end{array}\right.
 $$

So for a measurable function $u$ defined in $ M, T_{k}u$ is defined by $\left(T_{k} u\right)(x)=T_{k}(u(x)).$

we will use in its subsection Functional spaces :

\begin{enumerate}
\item[i)] $\tau_{l o c}^{1,1}(M)$ is the set of measurable functions $u : M \rightarrow \mathbb{R}$ such as for all $k>0$ the truncation function $T_{k}(u)$ in $W_{l o c}^{1,1}(M) .$
\item[ii)] for $p \in ] 1, \infty [, \tau_{l o c}^{1, p}(M)$ is the subset of $\tau_{l o c}^{1,1}(M)$ composed by functions $u$
such as $\left|\nabla\left(T_{k}(u)\right)\right| \in L_{l o c}^{p}(M)$ for all $k>0$.
\item[iii)] Of even, $\tau^{1, p}(M)$ is the subset of $\tau_{l o c}^{1,1}(M)$ composed of functions $u,$ such as, of
  more $\left|\nabla T_{k}(u)\right| \in L^{p}(M)$ for all $k>0$.
\item[iv)] Finally, $\tau_{0}^{1, p}(M)$ is the subset of $\tau^{1, p}(M),$ composed of functions that can be approximated by class functions $ C ^ {1} $ a compact support in $ M $ in the next sense : a function $u \in \tau^{1, p}(M)$ in  $\tau_{0}^{1, p}(M),$ if fopr all $k>0,$ it exists
a Sequence $\left(\phi_{n}\right) \subset C_{0}^{\infty}(M)$ such as
$$ 
\phi_{n} \rightarrow T_{k}(u) \quad \text { in } \quad L_{l o c}^{1}(M)
 $$
 $$ 
\nabla \phi_{n} \rightarrow \nabla T_{k}(u) \quad \text {in } \quad L^{p}(M)
 $$
\end{enumerate}
This space will play an important role in this work.\\

We have the next lemma giving some properties of the preceding spaces

\begin{lemma}
for all $p \in[1, \infty[,$ wa have
\begin{enumerate}
\item[1)] $W_{l o c}^{1, p}(M) \subset \tau_{l o c}^{1, p}(M) \quad$ et $\quad W_{0}^{1, p}(M) \subset \tau_{0}^{1, p}(M)$,
\item[2)] $\tau_{l o c}^{1, p}(M) \cap L_{l o c}^{\infty}(M)=W_{l o c}^{1, p}(M) \cap L_{l o c}^{\infty}(M)$,

\item[3)] $\nabla T_{k}(u)=\nabla u 1_{\{|u|<k\}},$

where $ {1} _ {A} $ denotes the characteristic function of a measurable set $A$.
\end{enumerate}
\end{lemma}

\textbf{Proof.}
\begin{enumerate}
\item[1)] we have
$$ 
\begin{aligned} u \in W_{l o c}^{1, p}(M) & \Rightarrow u \in W_{l o c}^{1,1}(M) \text { et } \nabla u \in L_{l o c}^{p}(M), \\ & \Rightarrow T_{k}(u) \in W_{l o c}^{1,1}(M) \quad \text { and } \nabla T_{k}(u) \in L_{l o c}^{p}(M) \quad \forall k>0 \\ & \Rightarrow u \in \tau_{l o c}^{1, p}(M) \end{aligned}
 $$
 so $W_{l o c}^{1, p}(M) \subset \tau_{l o c}^{1, p}(M) .$
For the second point, we have

\begin{center}
$u \in W_{0}^{1, p}(M) \Rightarrow u \in W^{1, p}(M) \quad$ and $\quad \exists\left\{\phi_{n}\right\} \subset C_{0}^{\infty}(M) \quad$ such as
\end{center}
$$\left\{ {\begin{array}{*{20}{c}}
{{\phi _n}\, \to u}&{{\rm{  }}}&{{\rm{in}}}&{{L^p}(M )}\\
{\nabla {\phi _n}}&{ \to \nabla u}&{{\rm{ in }}}&{{L^p}(M )}
\end{array}} \right.$$
\begin{center}
$ \Rightarrow u \in {\tau ^{1,p}}(M ){\rm{and}} \exists \left\{ {{\phi _n}} \right\} \subset C_0^\infty M )\quad such \,as$
\end{center}
$$ \Rightarrow \left\{ {\begin{array}{*{20}{c}}
{{\phi _n} \to {T_k}(u)}&{{\rm{ in }}}&{L_{loc}^1(M )}\\
{\nabla {\phi _n} \to \nabla {T_k}(u)}&{{\rm{ in }}}&{{L^p}(M )}
\end{array}} \right.\quad \forall k > 0,$$
$$ \Rightarrow u \in \tau _0^{1,p}(M ),$$
so $W_{0}^{1, p}(M) \subset \tau_{0}^{1, p}(M).$

\item[2)] as $$ 
\begin{aligned} u \in W_{l o c}^{1, p}(M) \cap L_{l o c}^{\infty}(M) & \Rightarrow u \in W_{l, o c}^{1, p}(M) \quad \text { and } \quad u \in L_{l o c}^{\infty}(M) \\ & \Rightarrow u \in \tau_{l o c}^{1, p}(M) \quad \text { and } \quad u \in L_{l o c}^{\infty}(M) \\ & \Rightarrow u \in \tau_{l o c}^{1, p}(M) \cap L_{l o c}^{\infty}(M) \end{aligned}
 $$
 so $\quad W_{l o c}^{1, p}(M) \cap L_{l o c}^{\infty}(M) \subset \tau_{l o c}^{1, p}(M) \cap L_{l o c}^{\infty}(\Omega).$
 
We also have
$$ 
\begin{aligned} u \in \tau_{l o c}^{1, p}(M) \cap L_{l o c}^{\infty}(M) & \Rightarrow \quad u \in \tau_{l o c}^{1, p}(M) \quad \text { and } \quad u \in L_{l o c}^{\infty}(M) \\ & \Rightarrow \quad T_{k}(u) \in W_{l o c}^{1,1}(M) \quad \text { and } \quad \nabla T_{k}(u) \in L_{l o c}^{p}(M) \quad \text { and } \quad u \in L_{l o c}^{\infty}(M) \\ & \Rightarrow \quad u \in L_{l o c}^{p}(M) \quad \text {and } \quad \nabla u \in L_{l o c}^{p}(M) \quad \text { and } \quad u \in L_{l o c}^{\infty}(M)\\ & 
\Rightarrow\quad u \in W_{l o c}^{1, p}(M) \quad \text{ and } \quad u \in L_{l o c}^{\infty}(M)\\ & 
\Rightarrow \quad u \in W_{l o c}^{1, p}(M) \cap L_{l o c}^{\infty}(M)
\end{aligned}
 $$
 so $\quad \tau_{l o c}^{1, p}(M) \cap L_{l o c}^{\infty}(M) \subset W_{l o c}^{1, p}(M) \cap L_{l o c}^{\infty}(M)$.
 
So $\quad \tau_{\text {loc}}^{1, p}(M) \cap L_{l o c}^{\infty}(M)=W_{l o c}^{1, p}(M) \cap L_{l o c}^{\infty}(M)$. 
\item[3)] We have  $$ 
T_{k}(u)=\left\{\begin{array}{ccc}{u} & {\text {if }} & {|u| \leq k} \\ {k \frac{u}{|u|}} & {\text {if }} & {|u|>k}\end{array}\right.
 $$

implies
$$\nabla {T_k}(u) = \left\{ \begin{array}{l}
\begin{array}{*{20}{c}}
{\nabla u}&\text{ if}&{\left| u \right| \le k}
\end{array}\\
\begin{array}{*{20}{c}}
{\nabla u}&\text{ if}&{\left| u \right| > k}
\end{array}
\end{array} \right.$$

so $\quad \nabla T_{k}(u)=\nabla u 1_{\{|u|<k\}}$.\hfill$\blacksquare$
\end{enumerate}

Note that if $u \in \tau_{l o c}^{1,1}(M),$ so $\nabla u$ is not defined even in the sense of distributions,
yet we have the next lemma that gives meaning to $\nabla u .$

\begin{lemma}$[1]$
Let $u \in \tau_{l o c}^{1,1}(M),$ it exists a function $v : M \rightarrow \mathbb{R}^{N}$ unique measurable such as 
$$ 
\nabla T_{k}(u)=v 1_{\{|u|<k\}} \quad a . e .$$
in others, $u \in W_{l o c}^{1,1}(M)$ 
if and only if
 $v \in L_{l o c}^{1}(M),$ so $v \equiv \nabla u$ in the usual weak sense.
\end{lemma}

\textbf{Proof. }

We have $\nabla T_{k}(u)=\nabla u 1_{\{|u|<<\}},$
so for all $u \in \tau_{l o c}^{1,1}(M)$ it exists a function $v : M \rightarrow \mathbb{R}^{N}$ measurable such as $v \equiv \nabla u$ a.e. and $v \in L_{l o c}^{1}(M)$.

$v$ is unique in the sense almost everywhere, because :

for all $k, \varepsilon>0,$ we have $T_{k}\left(T_{k+\varepsilon}(u)\right)=T_{k}(u) .$ Therefore, we get in $M_{k}=\{|u|<k\}$ legality $\nabla T_{k+\varepsilon}=\nabla T_{k}$ a.e.  hence the result, and so $v$ unique a.e.

It remains to show that $u \in W_{l o c}^{1,1}(M)$ if $v \in L_{l o c}^{1}(M) .$ Indeed, in this case $\nabla T_{k}(u) \rightarrow v$ in
$L_{l o c}^{1}(M),$ so we have to prove that $u \in L_{l o c}^{1}(M) .$ By contradiction, if $u \notin L_{l o c}^{1},$ there will be a closed ball $B \subset M$ such as
$$ 
t_{k}=\left\|T_{k}(u)\right\|_{L^{1}(B)} \rightarrow \infty \quad \text { when } \quad k \rightarrow \infty.
 $$
by
normalization, $v_{k}=\frac{T_{k}(u)}{t_{k}} .$ so $v_{k} \rightarrow 0$ a.e. $\left\|v_{k}\right\|_{L^{1}(B)}=1$ and $\left\|\nabla v_{k}\right\|_{L^{1}(B)} \rightarrow 0,$
contradiction with the compactness of the injection of $W^{1,1}(B)$ in $L^{1}(B) .$

\subsection{Solutions in the sense of entropy}
In this section we will develop the concept of the solution in the sense of entropy
which will allow us to study elliptic equations with second member in $L^{1}(M)$.

suppose that $f \in L^{1}(M)$ and consider the next equation:
\begin{equation} 
\left\{ {\begin{array}{*{20}{c}}
{ - {\Delta _p}u}& = &{f(x)}&{{\rm{ in }}}&M \\
u& = &0&{{\rm{ on }}}&{\partial M }
\end{array}} \right. 
\end{equation}
Let $u \in \tau_{0}^{1, p}(M)$ a solution of the equation (2) in $D^{\prime}(M),$ so for all $\phi \in C_{0}^{\infty}(M),$ we have 

$$ 
\int_{M}|\nabla u|^{p-2} \nabla u \nabla \phi d\sigma_{g}=\int_{M} f \phi d\sigma_{g} \quad \forall \phi \in C_{0}^{\infty}(M)
 $$

note that $f \in L^{1}(M),$ so by density and
if
we posit conditions of the type "Dirichlet homogeneous", so we can take $ T_{k}(u- \phi),k>0,$
as a test function in the previous equation we get

$$ 
\int_{M}|\nabla u|^{p-2} \nabla u \nabla T_{k}(u-\phi) d\sigma_{g}=\int_{M} T_{k}(u-\phi) f d\sigma_{g},
 $$
 
 so
 
 \begin{equation}
\int_{\{|u-\phi|<k\}}|\nabla u|^{ | p-2} \nabla u \nabla(u-\phi) d\sigma_{g}=\int_{M} T_{k}(u-\phi) f d\sigma_{g}.
\end{equation}

Note that each term in (3) is well defined $:$ as $\phi \in L^{\infty}(M),$ so

$$ 
\begin{aligned}|u-\phi|<k & \Rightarrow|u|-|\phi|<|u-\phi|<k \\ & \Rightarrow|u|<k+|\phi| \\ & \Rightarrow|u|<k+\|\phi\|_{\infty} \\ & \Rightarrow|u|<\overline{k} \end{aligned}
 $$ 
or $\overline{k}=k+\|\phi\|_{\infty}$.

as $|\nabla u|^{p-1} \in L^{1}(M),$ so

$$ 
\begin{aligned} \int_{\{|u-\phi|<k\}} &|\nabla u|^{p-2} \nabla u \nabla(u-\phi) d\sigma_{g} \\ &=\int_{\{|u-\phi|<k\}}|\nabla u|^{p} d\sigma_{g}-\int_{\{|u-\phi|<k\}}|\nabla u|^{p-2} \nabla u \nabla \phi d\sigma_{g} \\ & \leq \int_{\{|u|<\overline{k}\}}|\nabla u|^{p} d\sigma_{g}+\int_{\{|u-\phi|<k\}}|\nabla u|^{p-1}|\nabla \phi| d\sigma_{g}\\&
\leq \int_{M}\left|\nabla T_{\overline{k}}(u)\right|^{p}d\sigma_{g}+c_{1} \int_{\{|u-\phi|<k\}}|\nabla u|^{p} d\sigma_{g}+c_{2} \int_{\{|u-\phi|<k\}}|\nabla \phi|^{p} d\sigma_{g}\\&
\leq \quad c_{3} \int_{M}\left|\nabla T_{\overline{k}}(u)\right|^{p} d\sigma_{g}+c_{2} \int_{\{|u-\phi|<k\}}|\nabla \phi|^{p} d\sigma_{g}\\&
\leq C\left(\int_{M}\left|\nabla T_{\overline{k}}(u)\right|^{p} d\sigma_{g}+\int_{\{|u-\phi|<k\}}|\nabla \phi|^{p}d\sigma_{g}\right)
\end{aligned}
 $$
So

\begin{equation}
\int_{\{|u-\phi|<k\}}|\nabla u|^{p-2} \nabla u \nabla(u-\phi) d\sigma_{g} \leq C\left(\int_{M}\left|\nabla T_{\overline{k}}(u)\right|^{p} d\sigma_{g}+\int_{\{|u-\phi|<k\}}|\nabla \phi|^{p} d\sigma_{g}\right).
 \end{equation}

Since $T_{\overline{k}}(u) \in W_{0}^{1, p}(M)$ ie, $u \in \tau_{0}^{1, p}(M)$ and $\phi \in L^{\infty}(M) \cap W_{0}^{1, p}(M),$ the second
member in $(4)$ is bounded, so the first member of $ (3)$ is well defined.

We are in a position to give the next definition

\begin{definition}[Solution in the sense of entropy]
Let $f \in L^{1}(M),$ we say that $u \in\tau_{0}^{1, p}(M)$ is an entropy solution of the problem $(1)$ if $(3)$ is checked for each
$\phi \in L^{\infty}(M) \cap W_{0}^{1, p}(M)$ and for all $k>0$.
\end{definition}

Let's start by demonstrating some properties of the entropy solutions.

\begin{lemma}
Si $ u \in \tau_{0}^{1, p}(M)$ is an entropy solution of $(1)$ so for all $k>0$
$$ 
\frac{1}{k} \int_{\{|u|<k\}}|\nabla u|^{p} d\sigma_{g} \leq \int_{M}|f| d\sigma_{g}=\|f\|_{1}
 .$$
Therefore, we obtain the next estimate in $L^{p}(M)$
 \begin{equation}
\left\|\nabla T_{k}(u)\right\|_{p}^{p} \leq k\|f\|_{1}.
 \end{equation}
\end{lemma}

\textbf{Proof.}

As $\quad u \in \tau_{0}^{1, p}(M) \Rightarrow T_{k}(u) \in W_{0}^{1, p}(M) \Rightarrow T_{k}(u) \in L^{p}(M) . \quad$ If $\phi=0$ and grace at
$(3)$ we will have
$$ 
\int_{\{|u|<k\}}|\nabla u|^{p-2} \nabla u \nabla u d\sigma_{g}=\int T_{k}(u) f d\sigma_{g}=\int_{\{|u|<k\}} u f d\sigma_{g} \leq k \int_{\{|u|<k\}}|f| d\sigma_{g} \leq k \int_{M}|f| d\sigma_{g},
 $$
do $\displaystyle\int_{M}\left|\nabla T_{k}(u)\right|^{p} d\sigma_{g} \leq k\|f\|_{1}.$

\subsection{estimates}
Before demonstrating the existence of the entropy solution, we will prove some
preliminary estimates based on the estimate $ (5). $ These estimates will relate to $u$
and $ | \nabla u |$ in Marcinkiewicz spaces and we can consider them as keys to demonstrate compactness results in $L^{q}(M)$ spaces with $q$ suitably chosen. The first main result is the next lemma.

\begin{lemma}
Let $1<p<N$ and $ (M,g) $  a  Riemannian manifold of dimension $N$ Consider $u \in \tau_{0}^{1, p}(M)$ such as 
\begin{equation}
\frac{1}{k} \int_{\{|u|<k\}}|\nabla u|^{p} d\sigma_{g} \leq \alpha,
\end{equation}
for all $k>0 .$ So $u \in \mathcal{M}^{p_{1}}(M)$ with $p_{1}=\frac{N(p-1)}{N-p} .$ More precisely, there exists
$C=C(N, p)>0$ such as
\begin{equation}
\operatorname{meas}\{|u|>k\} \leq C \alpha^{\frac{N}{N-p}} k^{-p_{1}}.
\end{equation}
\end{lemma}

\textbf{Proof.}
Let $1<p<N$ and  $u \in \tau_{0}^{1, p}(M),$ so $T_{k}(u) \in W_{0}^{1, p}(M)$ for all $k>0,$ and according to the inequality of Sobolev we have

$$ 
\left\|T_{k}(u)\right\|_{p^{*}} \leq c(N, p)\left\|\nabla T_{k}(u)\right\|_{p} \quad \text { or } \quad p^{*}=\frac{N p}{N-p},
 $$
 
 because of $(6),$ we have  $\int_{M}\left|\nabla T_{k}(u)\right|^{p} d\sigma_{g} \leq k \alpha,$ ie. $\left\|\nabla T_{k}(u)\right\|_{p}^{p} \leq k \alpha$,
 
and consequently $\left\|\nabla T_{k}(u)\right\|_{p} \leq(k \alpha)^{\frac{1}{p}},$ so $\left\|T_{k}(u)\right\|_{p^{*}} \leq c(N, p)(k \alpha)^{\frac{1}{p}}$.

for $0<\varepsilon \leq k,$ We have $\{|u|>\varepsilon\}=\left\{\left|T_{k}(u)>\varepsilon\right|\right\},$ so

$$ 
\operatorname{meas}\{|u|>k\} \leq \varepsilon^{-p^{*}}\left\|T_{k}(u)\right\|_{p^{*}}^{p^{*}} \leq c_{1}(N, p)(k \alpha)^{\frac{p^{*}}{p}} \varepsilon^{-p^{*}} \leq c_{1}(N, p) \alpha^{\frac{N}{N-p}} k^{\frac{N}{N-p}} \varepsilon^{-\frac{N p}{N-p}}.
 $$

for $\varepsilon=k, \operatorname{meas}\{|u|>k\} \leq c_{1}(N, p) \alpha^{\frac{N}{N-p}} k^{-\frac{N(v-1)}{N-p}},$ we obtain
$$ 
\operatorname{meas}\{|u|>k\} \leq C k^{-p_{1}},
 $$
or $C=c_{1}(N, p) \alpha^{\frac{N}{N-p}} \quad$ and $p_{1}=\frac{N(p-1)}{N-p} .$ So it results than $\phi_{u}(k) \leq C k^{-p_{1}},$ and as a conclusion
it results than $u \in \mathcal{M}^{p_{1}}(M)$. \hfill$\blacksquare$

We now prove estimates on the gradient of $u$.

\begin{lemma}
Let $1<p<N$ and suppose that $u \in \tau_{0}^{1, p}(M)$ satisfied $(6)$ for all $k$.\\
 so for all $h>0$
$$ 
\operatorname{meas}\bigg\{|\nabla u|>h\bigg\} \leq C(N, p) \alpha^{\frac{N}{N-1}} h^{-p_{2}}, \quad p_{2}=\frac{N(p-1)}{N-1}
 .$$
\end{lemma}

\textbf{Proof}

for $k, \lambda>0,$ we pose
$$ 
\Phi(k, \lambda)=\operatorname{meas}\left\{|\nabla u|^{p}>\lambda,|u|>k\right\},
 $$
according to the Lemma $ (13) $ we have
\begin{equation}
\Phi(k, 0) \leq C(N, p) M^{\frac{N}{N-p}} h^{-p_{1}}
 \end{equation}

As the function $\lambda \mapsto \Phi(k, \lambda)$ is decreasing, we get for $k, \lambda>0$ and for 
$0 \leq s \leq \lambda, \Phi(0, \lambda) \leq \Phi(0, s),$ so

$$ 
\begin{aligned} \Phi(0, \lambda) \leq \Phi(0, s) & \Rightarrow \int_{0}^{\lambda} \Phi(0, \lambda) d s \leq \int_{0}^{\lambda} \Phi(0, s) d s \\ & \Rightarrow \Phi(0, \lambda) \leq \frac{1}{\lambda} \int_{0}^{\lambda} \Phi(0, s) d s \end{aligned}
 $$
 
 and
 
 $$ 
\begin{aligned} \frac{1}{\lambda} \int_{0}^{\lambda} \Phi(0, s) d s &=\frac{1}{\lambda} \int_{0}^{\lambda} \Phi(k, s) d s+\frac{1}{\lambda} \int_{0}^{\lambda}(\Phi(0, s)-\Phi(k, s)) d s \\ & \leq \frac{1}{\lambda} \int_{0}^{\lambda} \Phi(k, 0) d s+\frac{1}{\lambda} \int_{0}^{\lambda}(\Phi(0, s)-\Phi(k, s)) d s \\ & \leq \Phi(k, 0)+\frac{1}{\lambda} \int_{0}^{\lambda}(\Phi(0, s)-\Phi(k, s)) d s \end{aligned}
 $$
 
 so

\begin{equation}
\Phi(0, \lambda) \leq \frac{1}{\lambda} \int_{0}^{\lambda} \Phi(0, s) d s \leq \Phi(k, 0)+\int_{0}^{\lambda}(\Phi(0, s)-\Phi(k, s)) d s
\end{equation} 
 
 Note that
 
 $$ 
\Phi(0, s)-\Phi(k, s)=\operatorname{meas}\left\{|u|<k,|\nabla u|^{p}>s\right\}
 $$

as $(7)$, we will have

\begin{equation} 
\int_{0}^{\infty}(\Phi(0, s)-\Phi(k, s)) d s=\int_{\{|u|<k\}}|\nabla u|^{p} d\sigma_{g} \leq k \alpha.
 \end{equation}

Finally from $ (9) $ and using $ (8) $ and $ (10), $ we get to
$$ 
\Phi(0, \lambda) \leq \frac{\alpha k}{\lambda}+C(N, p) \alpha^{\frac{N}{N-p}} k^{-p_{1}}.
 $$

we pose $P(k)=\frac{\alpha k}{\lambda}+c \alpha^{\frac{N}{N-p}} k^{-p_{1}}$, so minimizing $P(k),$ of $k,$ we will have to solve the equation $P^{\prime}(k)=0,$ which implies that

$$ 
\frac{\alpha}{\lambda}-c p_{1} \alpha^{\frac{N}{N-p}} k^{-p_{1}-1}=0
 $$
 
 and so $k=\left(c \lambda p_{1} \alpha^{\frac{p}{N-p}}\right)^{\frac{1}{p_{1}+1}}.$
 
Consequently,
$$ 
\begin{aligned} \Phi(0, \lambda) & \leq k\left[\frac{\alpha}{\lambda}+c \alpha^{\frac{N}{N-p}} k^{-p_{1}-1}\right] \leq k\left[\frac{\alpha}{\lambda}+\frac{\alpha}{\lambda p_{1}} \alpha^{\frac{N}{N-p}} \alpha^{-\frac{N}{N-p}}\right] \\ & \leq k \frac{\alpha}{\lambda}\left[1+\frac{1}{p_{1}}\right] \leq \frac{\alpha}{\lambda}\left[1+\frac{1}{p_{1}}\right]\left(c \lambda p_{1} \alpha^{\frac{p}{N-p}}\right)^{\frac{1}{p_{1}+1}} \\&
\leq \frac{\alpha}{\lambda}\left[1+\frac{1}{p_{1}}\right]\left(c \lambda p_{1} \alpha^{\frac{p}{N-p}}\right)^{\frac{N-p}{p(N-1)}} \leq \frac{\alpha}{\lambda}\left[1+\frac{1}{p_{1}}\right]\left(c p_{1}\right)^{\frac{N-p}{p(N-1)}} \lambda^{-\frac{N-p}{p(N-1)}} \alpha^{\frac{p}{N-p}\frac{N-p}{p(N-1)}} \\&
 \leq\left[1+\frac{1}{p_{1}}\right]\left(c p_{1}\right)^{\frac{N-p}{p(N-1)}} \lambda^{-\frac{N(p-1)}{p(N-1)}} \alpha^{\frac{N}{N-1}}
\end{aligned}
 $$

so

$$ 
\Phi(0, \lambda) \leq C(N, p) \alpha^{\frac{N}{N-1}} \lambda^{-\frac{N(p-1)}{p(N-1)}}, \quad \text { with } \quad C(N, p)=\left[1+\frac{1}{p_{1}}\right]\left(c p_{1}\right)^{\frac{N-p}{p(N-1)}}
 $$

we pose $\lambda=h^{p},$ so

$$ 
\operatorname{meas}\left\{|\nabla u|^{p}>h^{p}\right\} \leq C(N, p) \alpha^{\frac{N}{N-1}} h^{-\frac{N(p-1)}{(N-1)}}
 $$
 
 and consequently
 
 $$ 
\operatorname{meas}\{|\nabla u|>h\} \leq C(N, p) \alpha^{\frac{N}{N-1}} h^{-p_{2}} \quad \text { with } \quad p_{2}=\frac{N(p-1)}{(N-1)}
 $$\hfill$\blacksquare$
 
 Hence the result.

\subsection{Existence of the entropy solution}

We are in a position to demonstrate the main result of this article, more precisely
we have the next theorem

\begin{theorem}
Let $1<p<N$ and Let $ (M,g) $  a compact Riemannian manifold , so it exists $u$ an entropy solution of the problem $(2)$ with $u \in \tau_{0}^{1, p}(M) .$ Furthermore
$$ 
u \in \mathcal{M}^{p_{1}}(M) \quad \text { and } \quad|\nabla u| \in \mathcal{M}^{p_{2}}(M),
 $$
or $p_{1}=\dfrac{N(p-1)}{N-p} \quad$ and $\quad p_{2}=\dfrac{N(p-1)}{N-1}.$

In the case $p>2-\dfrac{1}{N}$ the solution  $  
u \in
 W_{0}^{1, q}(M)$ for all $q<p_{2}$.
\end{theorem}

\textbf{Proof.}

The main idea of the demonstration is to proceed by approximation.

Step $1.$

As $f \in L^{1}(M)$ there is a sequence of functions $\left\{f_{n}\right\} \subset L^{\infty}(M)$ such as $f_{n} \longrightarrow f$ in  $L^{1}(M)$.

for $f_{n} \in L^{\infty}(M)$ it exists $u_{n} \in W_{0}^{1, p}(M),$  the unique weak solution of the problem

\begin{equation}
\left\{\begin{array}{ccc}{-\Delta_{p} u_{n}=f_{n}} & {\text { in }} & {M} \\ {u_{n}=0} & {\text { on }} & {\partial M}\end{array}\right.
\end{equation}

note that $T_{k}\left(u_{n}\right) \in L^{1}(M) \cap L^{\infty}(M)$ for all $k>0,$ so taking $T_{k}\left(u_{n}\right)$ as
test function in$(11)$ we get that

$$\int_{M}\left|\nabla u_{n}\right|^{p-2} \nabla u_{n} \nabla T_{k}\left(u_{n}\right) d\sigma_{g}=\int_{M} f_{n} T_{k}\left(u_{n}\right)d\sigma_{g}$$

so

$$ 
\int_{\left\{\left|u_{n}\right|<k\right\}}\left|\nabla u_{n}\right|^{p} d\sigma_{g} \leq k c
 $$

and so

$$ 
\frac{1}{k} \int_{\left\{\left|u_{n}\right|<k\right\}}\left|\nabla u_{n}\right|^{p} d\sigma_{g} \leq c
 $$

ie,

$$ 
\int_{M}\left|\nabla T_{k}\left(u_{n}\right)\right|^{p} d\sigma_{g} \leq k c
 $$

therefore it is concluded that $\left\{\nabla T_{k}\left(u_{n}\right)\right\}$ is bounded in $L^{p}(M)$ for all $k>0$ .
So it exists $w_{k}$ tel que $T_{k}\left(u_{n}\right) \rightarrow w_{k}$ weakly in $W_{0}^{1, p}(M)$ for each $k>0,$ and
$T_{k}\left(u_{n}\right) \rightarrow w_{k}$ a.e in $M .$ we pose $w_{k} \equiv T_{k}(u)$ in the set or $\left|w_{k}\right|<k,$ it's clear that $u$ is well defined because $T_{k+h}\left(u_{n}\right)=T_{k}\left(T_{h}\left(u_{n}\right)\right)$ and consequently $T_{k}\left(u_{n}\right) \rightarrow T_{k}(u)$
strongly in $L^{q}(M)$ for all $q<p^{*}$.

According to Lemmas 12 and $13,$ we have

\begin{center}
$u_{n} \in \mathcal{M}^{p_{1}}(M) \quad$ et $\quad\left|\nabla u_{n}\right| \in \mathcal{M}^{p_{2}}(M)$,
\end{center}

with $p_{1}=\frac{N(p-1)}{N-p}$ and $p_{2}=\frac{N(p-1)}{N-1}, \quad$ so
\begin{center}
$\left\|u_{n}\right\|_{\mathcal{M}^{p_{1}(M)}} \leq C \quad$ and $\quad\left\|\nabla u_{n}\right\|_{\mathcal{M}^{p_{2}}(M)} \leq \overline{C},$
\end{center}

and as $L^{q}(M) \subset \mathcal{M}^{q}(M) \subset L^{q-\varepsilon}(M)$ for all $q, \varepsilon>0,$ so

\begin{center}
$\left\|u_{n}\right\|_{L^{p_{1}-\varepsilon}(M)} \leq C \quad$ et $\quad\left\|\nabla u_{n}\right\|_{L^{p_{2}-\varepsilon}(M)} \leq \overline{C}.$
\end{center}

Note that if $p>2-\frac{1}{N},$ so $p_{2}>1$ and consequently $\left\{u_{n}\right\}$ will be bounded in
$W_{0}^{1, p_{2}-\varepsilon}(M)$ for all $\varepsilon>0$ with $p_{2}-\varepsilon \geq 1,$ so $u_{n} \rightarrow u$ weakly in $W_{0}^{1, p_{2}-\varepsilon}(M)$ .

So for all $\varphi \in W^{1, \infty}(M),$ we have

$$ 
\int_{M}\left|\nabla u_{n}\right|^{p-2} \nabla u_{n} \nabla \varphi d\sigma_{g}=\int_{M} f_{n} \varphi d\sigma_{g},
 $$
 
as $\left|\nabla u_{n}\right|^{p_{2}-\varepsilon} \equiv\left(\left|\nabla u_{n}\right|^{p-1}\right)^{\frac{N}{N-1}-\varepsilon}$ and $\quad p-1<p_{2},$ so for $\quad p-1<p_{2}-\varepsilon$
 $\left|\nabla u_{n}\right|^{p-1} \in L^{\frac{N}{N-1}-\varepsilon}(M)$.
 
as $|\nabla \varphi| \in L^{\infty}(M)$ and $f_{n} \rightarrow f \quad$ in $L^{1}(M),$ so going to the limit when
$n \rightarrow \infty,$ we find that

$$ 
\int_{M}|\nabla u|^{p-2} \nabla u \nabla \varphi d\sigma_{g}=\int_{M} f \varphi d\sigma_{g}
 $$

It's clear that $u_{n} \rightarrow u$ strongly in $L^{\overline{q}}(M)$ such as $1 \leq \overline{q}<\overline{p}^{*}$ with $\overline{p}^{*}=\frac{N \overline{p}}{N-\overline{p}}>1$
or $\overline{p}=p_{2}-\varepsilon$\\

\textbf{Step $2 .$ } To analyze the general case $1<p,$ we start by demonstrating that $u \in \tau_{0}^{1, p}(M) .$

We pose $\nabla T_{k}(u)=\nabla w_{k},$  is clear that $ \nabla T_{k}(u)$ is well defined because $w_{k} \in W_{0}^{1,p}(M),$ to go to the limit in $ k$ we will start by  show that $ \nabla u_{n}$ converges to $ \nabla u $ locally in measure. To prove it we show that $ \left \{\nabla u_{n} \right \}$ is a Cauchy sequence in measure.

Let $t$ and $\varepsilon>0,$ so

\begin{equation}
\begin{aligned}
\big\{\left|\nabla u_{n}-\nabla u_{m}\right|>t\big\} & \subset\big\{\left|\nabla u_{n}\right|>A\big\} \cup\big\{\left|\nabla u_{m}\right|>A\big\} \cup\big\{\left|u_{n}-u_{m}\right|>k\big\} \\ & \cup\big\{\left|u_{n}-u_{m}\right| \leq k,\left|\nabla u_{n}\right| \leq A,\left|\nabla u_{m}\right| \leq A,\left|\nabla u_{n}-\nabla u_{m}\right|>t\big\} .\end{aligned}
 \end{equation}

We choose $ A $ big enough as

$$ 
\operatorname{meas}\left\{\left|\nabla u_{n}\right|>A\right\} \leq \varepsilon \quad \text { for all } \quad n \in \mathbb{N},
 $$

(this is possible by Lemma 13).

To estimate the last term in $ (12), $ we use the next algebraic inequalities.

for all $\xi, \eta \in \mathbb{R}^{N},$ we have 

$$ 
\left\langle|\xi|^{p-2} \xi-|\eta|^{p-2} \eta, \xi-\eta\right\rangle \geq 0,
 $$

again if $ \ xi \ neq \ eta $ then
$$ 
\left\langle|\xi|^{p-2} \xi-|\eta|^{p-2} \eta, \xi-\eta\right\rangle> 0
 $$
 
and if $|\xi|<A,|\eta|<A$ and $|\xi-\eta|>t,$ so
it exists $ \mu>0 $ such that
$$\left\langle|\xi|^{p-2} \xi-|\eta|^{p-2} \eta, \xi-\eta\right\rangle \geq \mu .$$ 

Knowing that $-\Delta_{p} u_{n}=f_{n} \quad$ and $\quad-\Delta_{p} u_{m}=f_{m},$ so by subtracting and using
$T_{k}\left(u_{n}-u_{m}\right)$as a test function, we get

$$ 
\begin{aligned} \int_{\left\{\left|u_{n}-u_{m}\right| \leq k\right\}}\left\langle\left|\nabla u_{n}\right|^{p-2} \nabla u_{n}\right.&-\left|\nabla u_{m}\right|^{p-2} \nabla u_{m}, \nabla u_{n}-\nabla u_{m} \rangle d\sigma_{g} \\ &=\int_{M}\left(f_{n}-f_{m}\right) T_{k}\left(u_{n}-u_{m}\right) d\sigma_{g} \leq 2 c k. \end{aligned}
 $$

According to the Lemma $ 12, $ we have

$$ 
\begin{aligned} & \quad \operatorname{meas}\left\{\left|u_{n}-u_{m}\right| \leq k,\left|\nabla u_{n}\right| \leq A,\left|\nabla u_{m}\right| \leq A,\left|\nabla u_{n}-\nabla u_{m}\right|>t\right\} \\ \leq & \quad \operatorname{meas}\left\{\left|u_{n}-u_{m}\right| \leq k, \quad\left(\left|\nabla u_{n}\right|^{p-2} u_{n}-\left|\nabla u_{m}\right|^{p-2} u_{m}\right) \cdot\left(\nabla u_{n}-\nabla u_{m}\right) \geq \mu\right\} \\ \leq & \quad \frac{1}{\mu} \int_{\left\{\left|u_{n}-u_{m}\right| \leq k\right\}}\left\langle\left|\nabla u_{n}\right|^{p-2} u_{n}-\left|\nabla u_{m}\right|^{p-2} u_{m}, \nabla u_{n}-\nabla u_{m}\right\rangle d\sigma_{g} \\ \leq & \quad \frac{1}{\mu} 2 c k \leq \varepsilon, \end{aligned}
 $$
if $k$ is small enough, as $k \leq \dfrac{\mu \varepsilon}{2 c}$.

So we fix $ A $ and $ k, $ if $ n_ {0} $ big enough, we have to $n, m \geq n_{0}, \operatorname{mes}\big\{\left|u_{n}-u_{m}\right|> k \big\} \leq  \varepsilon,$ and so $$\operatorname{meas}\left\{\left|\nabla u_{n}-\nabla u_{m}\right|>k\right\} \leq 2 \varepsilon.$$

So $\left\{\nabla u_{n}\right\}$ converges locally to a $ v $ function and as a consequence
a.e. in $ M. $ Since $ \left \{\nabla T_ {k} \left (u_ {n} \right) \right \} $ is bounded in $ L ^ {p} (M) $ for all $ k> 0 $ and $ \nabla T_ {k} \left (u_ {n} \right) \rightharpoonup $
$ \nabla T_ {k} (u) $ weakly in $ L ^ {p} (M), $ we deduce that $ v = \nabla u $ a.e. Note that in general
$ v \notin \left (L ^ {1} (M) \right) ^ {N}. $ It is clear that if $ p> 2- \frac {1} {N}, $ then $ v \notin \left (L ^ {1} (M) \right) ^ {N}, $ and so $ u \in W_ {0} ^ {1,1} (M) $
and from Lemma 10 we deduce $ \nabla u = v \quad a. e.  $

And consequently $ u \in \tau_ {0} ^ {1,1} (M) $.

To see that $ u \in \tau_ {0} ^ {1, p} (M), $ we consider $ \phi_ {n} \in C_ {0} ^ {\infty} (M) $ such that
$$ 
\left\|\nabla \phi_{n}-\nabla T_{k}\left(u_{n}\right)\right\|_{L^{p}(\Omega)} \leq \frac{1}{n} \quad \text { et } \quad\left\|\phi_{n}-T_{k}\left(u_{n}\right)\right\|_{L^{p^{*}}(\Omega)} \leq \frac{1}{n}.
 $$

We  have then
$$ 
\nabla \phi_{n} \longrightarrow \nabla T_{k}(u) \quad \text { fortement dans } \quad L^{p}(M)
 $$

and

$$ 
\phi_{n} \longrightarrow T_{k}(u) \quad \text { fortement in } \quad L_{l o c}^{q}(M) \quad \text { for } \quad q<p^{*}.
 $$

As a conclusion we get that $ \phi_ {n} $ converges strongly to $ T_ {k} (u) $ and consequently
$u \in \tau_{0}^{1, p}(M) .$

\textbf{Step $ 3. $} In this step we will demonstrate the strong convergence of truncations in
$ W_ {0} ^ {1, p} (M), $ ie for $ k> 0 $ fixed on a $ T_ {k} \left (u_ {n} \right) \rightarrow T_ {k} (u) $ strongly in $ W_ {0} ^ {1, p} (M). $

Note that $ T_ {k} \left (u_ {n} \right) \rightarrow T_ {k} (u) $ weakly in $ W_ {0} ^ {1, p} (M) $ for all $ k> 0 $.

Let $ k, h> 0 $ such that $ h> k> 0, $
we assume

$$ w_ {n} = T_ {2 k} \left (u_ {n} -T_ {h} \left (u_ {n} \right) T_ {k} \left (u_ {n} \right) -T_ {k} (u) \right) $$
Taking $ w_ {n} $ as a test function in $ (11), $ it results

$$ 
\int_{M}\left|\nabla u_{n}\right|^{p-2} \nabla u_{n} \nabla w_{n} d\sigma_{g}=\int_{M} f_{n} w_{n} d\sigma_{g},
 $$

we pose $I=\int_{M}\left|\nabla u_{n}\right|^{p-2} \nabla u_{n} \nabla w_{n} d x,$ when $k \rightarrow \infty$ and $h \rightarrow \infty$ we have $\int_{M} f_{n} w_{n} d x \rightarrow 0,$
so $\int_{M}\left|\nabla u_{n}\right|^{p-2} \nabla u_{n} \nabla w_{n} d x \rightarrow 0.$

We pose $\alpha=4 k+h . \quad$ if $\left|u_{n}\right|>\alpha, \nabla w_{n}=0 .$ So

$$ 
\begin{aligned} I &=\int_{\left\{\left|u_{n}\right|<\alpha\right\}}\left|\nabla u_{n}\right|^{p-2} \nabla u_{n} \nabla w_{n} d\sigma_{g} \\ &=\int_{\left\{\left|u_{n}\right|<k\right\}}\left|\nabla u_{n}\right|^{p-2} \nabla u_{n} \nabla w_{n} d\sigma_{g}+\int_{\left\{k<\left|u_{n}\right|<\alpha\right\}}\left|\nabla u_{n}\right|^{p-2} \nabla u_{n} \nabla w_{n} d\sigma_{g} \\ &=\int_{M}\left|\nabla T_{k}\left(u_{n}\right)\right|^{p-2} \nabla T_{k}\left(u_{n}\right) \nabla\left(T_{k}\left(u_{n}\right)-T_{k}(u)\right) d\sigma_{g} \\ &+\int_{\left\{k<\left|u_{n}\right|<\alpha\right\}}\left|\nabla u_{n}\right|^{p-2} \nabla u_{n} \nabla w_{n} d\sigma_{g}. \end{aligned}
 $$

note that
$1-$ if $k<\left|u_{n}\right| \leq h$ so $\nabla w_{n}=\nabla T_{k}(u)$

and

$2-$ if $h<\left|u_{n}\right|<\alpha$ so $\nabla w_{n}=\nabla T_{k}(u)$,

so

$$ 
\begin{aligned} \int_{\left\{k<\left|u_{n}\right|<\alpha\right\}}\left|\nabla u_{n}\right|^{p-2} \nabla u_{n} \nabla w_{n} d\sigma_{g} &=\int_{\left\{\left|u_{n}\right|>k\right\}}\left|\nabla T_{\alpha}\left(u_{n}\right)\right|^{p-2} \nabla T_{\alpha}\left(u_{n}\right) \nabla T_{k}(u) d\sigma_{g} \\ & \geq-\int_{\left\{\left|u_{n}\right|>k\right\}}\left|\nabla T_{\alpha}\left(u_{n}\right)\right|^{p-1}\left|\nabla T_{k}(u)\right| d\sigma_{g} \end{aligned}
 $$
 
 we obtain
 
 $$ 
\begin{aligned} I & \geq \int_{M}\left|\nabla T_{k}\left(u_{n}\right)\right|^{p-2} \nabla T_{k}\left(u_{n}\right) \nabla\left(T_{k}\left(u_{n}\right)-T_{k}(u)\right) d\sigma_{g} \\ &-\int_{\left\{\left|u_{n}\right|>k\right\}}\left|\nabla T_{\alpha}\left(u_{n}\right)\right|^{p-1}\left|\nabla T_{k}(u)\right| d\sigma_{g}. \end{aligned}
 $$

As $\left\{\left|\nabla T_{\alpha}\left(u_{n}\right)\right|^{p-1}\right\}$ is bounded in $L^{\frac{p}{p-1}}(M),\left|\nabla T_{k}(u)\right| 1_{\left\{\left|u_{n}\right|>k\right\}}$ is bounded in
$L^{p}(M)$ and $\left|\nabla T_{k}(u)\right| \mathbb{1}_{\left\{\left|u_{n}\right|>k\right\}} \rightarrow 0$ strongly in $L^{p}(M)$ when $k \rightarrow \infty,$ we get that

$$ 
\int_{\left\{\left|u_{n}\right|>k\right\}}\left|\nabla T_{\alpha}\left(u_{n}\right)\right|^{p-1}\left|\nabla T_{k}(u)\right| d\sigma_{g} \rightarrow 0 \quad \text { when } \quad k \rightarrow \infty .
 $$

So

$$ 
\begin{aligned} J &=\int_{M}\left|\nabla T_{k}\left(u_{n}\right)\right|^{p-2} \nabla T_{k}\left(u_{n}\right) \nabla\left(T_{k}\left(u_{n}\right)-T_{k}(u)\right) d x \\ &=\int_{M}\left(\left|\nabla T_{k}\left(u_{n}\right)\right|^{p-2} \nabla T_{k}\left(u_{n}\right)-\left|\nabla T_{k}(u)\right|^{p-2} \nabla T_{k}(u)\right)\left(\nabla T_{k}\left(u_{n}\right)-\nabla T_{k}(u)\right) d\sigma_{g} \\ &+\int_{M}\left|\nabla T_{k}(u)\right|^{p-2} \nabla T_{k}(u)\left(\nabla T_{k}\left(u_{n}\right)-\nabla T_{k}(u)\right) d\sigma_{g} \end{aligned}
 $$
 
as
$$ 
\int_{M}\left|\nabla T_{k}(u)\right|^{p-2} \nabla T_{k}(u)\left(\nabla T_{k}\left(u_{n}\right)-\nabla T_{k}(u)\right) d\sigma_{g} \rightarrow 0 \quad \text { when } \quad n \rightarrow \infty ,
 $$

so
$$ 
J=\int_{M}\left(\left|\nabla T_{k}\left(u_{n}\right)\right|^{p-2} \nabla T_{k}\left(u_{n}\right)-\left|\nabla T_{k}(u)\right|^{p-2} \nabla T_{k}(u)\right)\left(\nabla T_{k}\left(u_{n}\right)-\nabla T_{k}(u)\right) d\sigma_{g}+o(1).
 $$
So
$$ 
\begin{aligned} I & \geq \int_{M}\left(\left|\nabla T_{k}\left(u_{n}\right)\right|^{p-2} \nabla T_{k}\left(u_{n}\right)-\left|\nabla T_{k}(u)\right|^{p-2} \nabla T_{k}(u)\right)\left(\nabla T_{k}\left(u_{n}\right)-\nabla T_{k}(u)\right) d\sigma_{g}+o(1) \\ & \geq c \int_{M}\left|\nabla T_{k}\left(u_{n}\right)-\nabla T_{k}(u)\right|^{p} d\sigma_{g}+o(1) \quad \text {if } \quad p \geq 2 \end{aligned}
 $$

and
$$ 
I \geq C(p) \int_{M} \frac{\left|\nabla T_{k}\left(u_{n}\right)-\nabla T_{k}(u)\right|^{2}}{\left(\left|\nabla T_{k}(u)\right|+\left|\nabla T_{k}\left(u_{n}\right)\right|\right)^{2-p}} d\sigma_{g}+o(1) \quad \text { if } \quad p<2.
 $$

So
$$ 
\int_{M}\left|\nabla T_{k}\left(u_{n}\right)-\nabla T_{k}(u)\right|^{p} d\sigma_{g} \leq o(1)+\int_{M} f_{n} w_{n} d\sigma_{g} \quad \text { if } \quad p \geq 2
 $$

and
$$ 
C(p) \int_{M} \frac{\left|\nabla T_{k}\left(u_{n}\right)-\nabla T_{k}(u)\right|^{2}}{\left(\left|\nabla T_{k}(u)\right|+\left|\nabla T_{k}\left(u_{n}\right)\right|\right)^{2-p}} d\sigma_{g} \leq o(1)+\int_{M} f_{n} w_{n} \quad \text { if } \quad p<2.
 $$

Consequently
$$ 
\int_{M}\left|\nabla T_{k}\left(u_{n}\right)-\nabla T_{k}(u)\right|^{p} d\sigma_{g} \rightarrow 0 \quad \text { if } \quad p \geq 2, $$

and 
$$ 
\int_{M} \frac{\left|\nabla T_{k}\left(u_{n}\right)-\nabla T_{k}(u)\right|^{2}}{\left(\left|\nabla T_{k}(u)\right|+\left|\nabla T_{k}\left(u_{n}\right)\right|\right)^{2-p}} d\sigma_{g} \rightarrow 0 \quad \text { if } \quad p<2.
 $$

For the second case we have
$$ 
\begin{array}{c}{\displaystyle\int_{M}\left|\nabla T_{k}\left(u_{n}\right)-\nabla T_{k}(u)\right|^{p} d\sigma_{g}=\displaystyle\int_{M} \frac{\left|\nabla T_{k}\left(u_{n}\right)-\nabla T_{k}(u)\right|^{p}}{\left(|\nabla T_{k}(u)|+| \nabla T_{k}\left(u_{n}\right) |\right)^{\frac{p(2-p)}{2}}}\left(\left|\nabla T_{k}(u)\right|+\left|\nabla T_{k}\left(u_{n}\right)\right|\right)^{\frac{p(2-p)}{2}} d\sigma_{g}} \\ {\quad \leq\left(\displaystyle\int_{M} \frac{\left|\nabla T_{k}\left(u_{n}\right)-\nabla T_{k}(u)\right|^{2}}{\left(\left|\nabla T_{k}(u)\right|+\left|\nabla T_{k}\left(u_{n}\right)\right|\right)^{2-p}} d\sigma_{g}\right)^{\frac{p}{2}}\left(\displaystyle\int_{M}\left(\left|\nabla T_{k}\left(u_{n}\right)\right|^{p}+\left|\nabla T_{k}(u)\right|^{p}\right) d\sigma_{g}\right)^{\frac{2-p}{2}}}.\end{array}
 $$

So for $p < 2$,
$$ 
\int_{M}\left|\nabla T_{k}\left(u_{n}\right)-\nabla T_{k}(u)\right|^{p} d\sigma_{g} \rightarrow 0 \text { for } n \rightarrow \infty
 $$

As a conclusion we obtain that $T_{k}\left(u_{n}\right) \rightarrow T_{k}(u)$ strongly in $W_{0}^{1, p}(M)$ for all $k>0 .$

\textbf{Step $ 4. $} To complete the proof it remains to show that $ u $ is an entropy solution.

recall that 
$$ 
\left\{\begin{array}{ccc}{-\Delta_{p} u_{n}=f_{n}} & {\text { in }} & {M} \\ {u_{n}=0} & {\text { on}} & {\partial M}\end{array}\right.
 $$

Let $v \in L^{\infty}(M) \cap W_{0}^{1, p}(M),$ for all $k$ fixed $>0,$ we have 
$$\int_{M}\left|\nabla u_{n}\right|^{p-2} \nabla u_{n} \nabla T_{k}\left(u_{n}-v\right) d\sigma_{g}=\int_{M} f_{n} T_{k}\left(u_{n}-v\right) d\sigma_{g}$$

As $u_{n} \longrightarrow u \quad a . e .$ in $M,$ and $f_{n} \rightarrow f \quad$ in $L^{1}(M) .$

So $\int_{M} f_{n} T_{k}\left(u_{n}-v\right) d\sigma_{g} \rightarrow \int_{M} f T_{k}(u-v) d\sigma_{g}$ for $n \rightarrow \infty$.

As $v \in L^{\infty}(M) \cap W_{0}^{1, p}(M),$ so it exists a positive constant $ c> 0 $ such that

$$ 
\int_{M}\left|\nabla u_{n}\right|^{p-2} \nabla u_{n} \nabla T_{k}\left(u_{n}-v\right) d\sigma_{g}=\int_{\left\{\left|u_{n}\right| \leq c\right\}}\left|\nabla T_{c}\left(u_{n}\right)\right|^{p-2} \nabla T_{c}\left(u_{n}\right) \nabla T_{k}\left(u_{n}-v\right) d\sigma_{g}
 $$
 
 Note that it is sufficient to take $c \geq k+\|v\|_{\infty} .$ As $T_{k}\left(u_{n}\right) \rightarrow T_{k}(u)$ strongly in
$W_{0}^{1, p}(M),$ so we conclude that

$$ 
\int_{\left\{\left|u_{n}\right| \leq c\right\}}\left|\nabla T_{c}\left(u_{n}\right)\right|^{p-2} \nabla T_{c}\left(u_{n}\right) \nabla T_{k}\left(u_{n}-v\right) d\sigma_{g} \rightarrow \int_{\{|u| \leq c\}}\left|\nabla T_{c}(u)\right|^{p-2} \nabla T_{c}(u) \nabla T_{k}(u-v) d\sigma_{g}
 $$

for $n \rightarrow \infty .$ Consequently and for $n \rightarrow \infty$ we get that
$$\int_{M}|\nabla u|^{p-2} \nabla u \nabla T_{k}(u-v) d\sigma_{g}=\int_{M} f T_{k}(u-v)d\sigma_{g}$$
So $u$ is an entropy solution of the problem $(1,2).$ \hfill$\blacksquare$

\subsection{Uniqueness of the solution in the sense of entropy}

We deal here with the question of the uniqueness of entropy solutions $ u \in \tau_ {0} ^ {1, p} (M) $ for the problem $ (2), $ note that $ u $ checks $ (3) $ for all $ \phi \in L^{\infty} (M) \cap W_ {0}^{1, p}(M) $
and for all $ k>0 $

The main result of this section is the next theorem

\begin{theorem}
Let $u_{1}$ et $u_{2}$ dtwo functions in $\tau_{0}^{1, p}(M),$ such as $u_{1}$ and $u_{2}$ are entropy solutions to the problem
$$-\Delta_{p} u=f(x)$$
so $u_1 = u_2$.
\end{theorem}

\textbf{Proof.}

Note that $ f \in L ^ {1} (M), $
substitute
in the relation $ (3) $ with test functions
$ T_ {h} \left (u_ {1} \right) $ and $ T_ {h} \left (u_ {2} \right) $ and by addition gets that

$$ 
\int_{\left\{\left|u_{1}-T_{h}\left(u_{2}\right)\right|<k\right\}}\left|\nabla u_{1}\right|^{p-2} \nabla u_{1} \nabla\left(u_{1}-T_{h}\left(u_{2}\right)\right) d\sigma_{g}=\int_{M} f T_{k}\left(u_{1}-T_{h}\left(u_{2}\right)\right) d\sigma_{g},$$

$$\int_{\left\{\left|u_{2}-T_{h}\left(u_{1}\right)\right|<k\right\}}\left|\nabla u_{2}\right|^{p-2} \nabla u_{2} \nabla\left(u_{2}-T_{h}\left(u_{1}\right)\right) d\sigma_{g}=\int_{M} f T_{k}\left(u_{2}-T_{h}\left(u_{1}\right)\right)d\sigma_{g}.
 $$

By combining the two results we get\\

$\displaystyle\int_{\left\{\left|u_{1}-T_{h}\left(u_{2}\right)\right|<k\right\}}\left|\nabla u_{1}\right|^{p-2} \nabla u_{1} \nabla\left(u_{1}-T_{h}\left(u_{2}\right)\right) d\sigma_{g}$

\begin{equation}
\begin{aligned} &+\int_{\left\{\left|u_{2}-T_{h}\left(u_{1}\right)\right|<k\right\}}\left|\nabla u_{2}\right|^{p-2} \nabla u_{2} \nabla\left(u_{2}-T_{h}\left(u_{1}\right)\right) d\sigma_{g} \\ &=\int_{M} f\left(T_{k}\left(u_{1}-T_{h}\left(u_{2}\right)\right)+T_{k}\left(u_{2}-T_{h}\left(u_{1}\right)\right)\right) d\sigma_{g}. \end{aligned}
 \end{equation}

The conclusion $ u_ {1} = u_ {2} $
will be reached after going to the limit $ h \ rightarrow \ infty $ in this
formula.
Let

$$ 
\begin{aligned} I &=\int_{\left\{\left|u_{1}-T_{h}\left(u_{2}\right)\right|<k\right\}}\left|\nabla u_{1}\right|^{p-2} \nabla u_{1} \nabla\left(u_{1}-T_{h}\left(u_{2}\right)\right) d\sigma_{g} \\ &+\int_{\left\{\left|u_{2}-T_{h}\left(u_{1}\right)\right|<k\right\}}\left|\nabla u_{2}\right|^{p-2} \nabla u_{2} \nabla\left(u_{2}-T_{h}\left(u_{1}\right)\right)d\sigma_{g}. \end{aligned}
 $$
 
 we pose
 
 $$ 
A_{0}=\big\{x \in M :\left|u_{1}-u_{2}\right|<k,\left|u_{1}\right|<h,\left|u_{2}\right|<h\big\}.
 $$

In $ A_ {0} $ the first member of $ (13) $ is reduced to the next term

$$I_{0}=\int_{A_{0}}\left(\left|\nabla u_{1}\right|^{p-2} \nabla u_{1}-\left|\nabla u_{2}\right|^{p-2} \nabla u_{2}\right)\left(\nabla u_{1}-\nabla u_{2}\right) d\sigma_{g}.$$

Let$$ 
A_{1}=\left\{x \in M :\left|u_{1}-T_{h}\left(u_{2}\right)\right|<k,\left|u_{2}\right| \geq h\right\},
 $$

so

$$ 
\int_{A_{1}}\left|\nabla u_{1}\right|^{p-2} \nabla u_{1} \nabla\left(u_{1}-T_{h}\left(u_{2}\right)\right) d\sigma_{g}=\int_{A_{1}}\left|\nabla u_{1}\right|^{p}d\sigma_{g} \geq 0,
 $$

and on set

$$ 
A_{2}=\left\{x \in M :\left|u_{1}-T_{h}\left(u_{2}\right)\right|<k,\left|u_{2}\right|<h,\left|u_{1}\right| \geq h\right\},
 $$

we are getting

$$ 
\begin{aligned} \int_{A_{2}}\left|\nabla u_{1}\right|^{p-2} \nabla u_{1} \nabla\left(u_{1}-T_{h}\left(u_{2}\right)\right) d\sigma_{g} &=\int_{A_{2}}\left|\nabla u_{1}\right|^{p-2} \nabla u_{1}\left(\nabla u_{1}-\nabla u_{2}\right) d\sigma_{g} \\ & \geq-\int_{A_{2}}\left|\nabla u_{1}\right|^{p-2} \nabla u_{1} \nabla u_{2} d\sigma_{g}. \end{aligned}
 $$

In the same way, we can define all $A_{1}^{\prime}$ and $A_{2}^{\prime}$ as
$$A_{1}^{\prime}=\left\{x \in M :\left|u_{2}-T_{h}\left(u_{1}\right)\right|<k,\left|u_{1}\right| \geq h\right\},$$

and

$$ 
A_{2}^{\prime}=\left\{x \in M :\left|u_{2}-T_{h}\left(u_{1}\right)\right|<k,\left|u_{1}\right|<h,\left|u_{2}\right| \geq h\right\}
 $$

Then the second term of $ (13) $ can be written as a sum of 

$$ 
\int_{A_{1}^{\prime}}\left|\nabla u_{2}\right|^{p-2} \nabla u_{2}\left(\nabla u_{2}-\nabla T_{h}\left(u_{1}\right)\right) d\sigma_{g}=\int_{A_{1}^{\prime}}\left|\nabla u_{2}\right|^{p} d\sigma_{g} \geq 0
 $$
 
 and
 
 $$ 
\begin{aligned} \int_{A_{2}^{\prime}}\left|\nabla u_{2}\right|^{p-2} \nabla u_{2}\left(\nabla u_{2}-\nabla T_{h}\left(u_{1}\right)\right)d\sigma_{g} &=\int_{A_{2}^{\prime}}\left|\nabla u_{2}\right|^{p-2} \nabla u_{2}\left(\nabla u_{2}-\nabla u_{1}\right) d\sigma_{g} \\ & \geq-\int_{A_{2}^{\prime}}\left|\nabla u_{2}\right|^{p-2} \nabla u_{2} \nabla u_{1} d\sigma_{g} \end{aligned}
 $$
 
 Therefore we conclude that
 
 $$ 
\begin{aligned} I & \geq I_{0}+\int_{A_{1}}\left|\nabla u_{1}\right|^{p} d\sigma_{g}-\int_{A_{2}}\left|\nabla u_{1}\right|^{p-2} \nabla u_{1} \nabla u_{2} d\sigma_{g} \\ &+\int_{A_{1}^{\prime}}\left|\nabla u_{2}\right|^{p} d\sigma_{g}-\int_{A_{2}^{\prime}}\left|\nabla u_{2}\right|^{p-2} \nabla u_{2} \nabla u_{1} d\sigma_{g} \\ & \geq I_{0}-\left(\int_{A_{2}}\left|\nabla u_{1}\right|^{p-2} \nabla u_{1} \nabla u_{2} d\sigma_{g}+\int_{A_{2}^{\prime}}\left|\nabla u_{2}\right|^{p-2} \nabla u_{2} \nabla u_{1} d\sigma_{g}\right) \\ & \geq I_{0}-I_{3} \end{aligned}
 $$
 
or
 
 $$ 
I_{3}=\int_{A_{2}}\left|\nabla u_{1}\right|^{p-2} \nabla u_{1} \nabla u_{2} d\sigma_{g}+\int_{A_{2}^{\prime}}\left|\nabla u_{2}\right|^{p-2} \nabla u_{2} \nabla u_{1} d\sigma_{g}
 $$
 
The first term of $ I_ {3} $ can be estimated by

$$ 
\begin{aligned} \int_{A_{2}}\left|\nabla u_{1}\right|^{p-2} \nabla u_{1} \nabla u_{2} d\sigma_{g} & \leq \int_{A_{2}}\left|\nabla u_{1}\right|^{p-1}\left|\nabla u_{2}\right| d\sigma_{g} \\ & \leq\left(\int_{A_{2}}\left|\nabla u_{1}\right|^{p} d\sigma_{g}\right)^{\frac{p-1}{p}}\left(\int_{A_{2}}\left|\nabla u_{2}\right|^{p} d\sigma_{g}\right)^{\frac{1}{p}} \\ & \leq
\left\|\nabla u_{1}\right\|_{L^{p}\left(\left\{h \leq\left|u_{1}\right| \leq h+k\right\}\right)}^{p-1}\left\|\nabla u_{2}\right\|_{L^{p}\left(\left\{h-k \leq\left|u_{2}\right| \leq h\right\}\right)}. 
\end{aligned}
 $$

as $\left\|\nabla u_{1}\right\|_{L^{p}\left(\left\{h \leq\left|u_{1}\right| \leq h+k\right\}\right)}^{p-1}\left\|\nabla u_{2}\right\|_{L^{p}\left(\left\{h-k \leq\left|u_{2}\right| \leq h\right\}\right)} \rightarrow 0$ when $h \rightarrow \infty$for all $k>0,$ it results
than $\displaystyle\int_{A_{2}}\left|\nabla u_{1}\right|^{p-2} \nabla u_{1} \nabla u_{2} d\sigma_{g}$ converges to $0$ when $h \rightarrow \infty$for tout $k>0$.

In the same way we obtain the same conclusion for the second term of $I_{3}$ .

So we conclude that $ I_ {3} $ tends to 0 when $h \rightarrow \infty$ .\\

Regarding the second member of $ (13), $ knowing that

\begin{center}
$T_{k}\left(u_{1}-T_{h}\left(u_{2}\right)\right)+T_{k}\left(u_{2}-T_{h}\left(u_{1}\right)\right) \rightarrow 0$ a.e. in $in$ for $h \rightarrow \infty$
\end{center}
$$\left|T_{k}\left(u_{1}-T_{h}\left(u_{2}\right)\right)+T_{k}\left(u_{2}-T_{h}\left(u_{1}\right)\right)\right| \leq 2 k$$

and that $ f \ in L ^ {1} (M), $ so using the dominated Convergence Theorem we get that

$$ 
\int_{M} f\left(T_{k}\left(u_{1}-T_{h}\left(u_{2}\right)\right)+T_{k}\left(u_{2}-T_{h}\left(u_{1}\right)\right)\right) d\sigma_{g}\rightarrow 0 \text { quand } h \rightarrow \infty \text { for all } k>0.
 $$

Combining previous estimates
it results
than

$$ 
\int_{A_{0}(h, k)}\left(\left|\nabla u_{1}\right|^{p-2} \nabla u_{1}-\left|\nabla u_{2}\right|^{p-2} \nabla u_{2}\right)\left(\nabla u_{1}-\nabla u_{2}\right) d\sigma_{g} \leq \varepsilon(h),
 $$
or $\varepsilon(h) \rightarrow 0$ when $h \rightarrow \infty$ for all $k$ fixed $>0 .$ Since $A_{0}(h, k)$ converges to
$$\left\{x \in M :\left|u_{1}-u_{2}\right|<k\right\},$$

we conclude that
$$ 
\int_{\left\{\left|u_{1}-u_{2}\right|<k\right\}}\left(\left|\nabla u_{1}\right|^{p-2} \nabla u_{1}-\left|\nabla u_{2}\right|^{p-2} \nabla u_{2}\right)\left(\nabla u_{1}-\nabla u_{2}\right) d\sigma_{g} \leq 0.
 $$
As

$$ 
\lambda\left\|\nabla u_{1}-\nabla u_{2}\right\|_{L^{p}\left(\left\{\left|u_{1}-u_{2}\right|<k\right\}\right)}^{p} \leq \int_{\left\{\left|u_{1}-u_{2}\right|<k\right\}}\left(\left|\nabla u_{1}\right|^{p-2} \nabla u_{1}-\left|\nabla u_{2}\right|^{p-2} \nabla u_{2}\right)\left(\nabla u_{1}-\nabla u_{2}\right) d\sigma_{g}
 $$

if $p > 2$ and

$$ 
\int_{\left\{\left|u_{1}-u_{2}\right|<k\right\}} \frac{\left|\nabla u_{1}-\nabla u_{2}\right|^{2}}{\left(\left|\nabla u_{1}\right|+\left|\nabla u_{2}\right|\right)^{2-p}} d\sigma_{g} \leq $$

$$
\int_{\left\{\left|u_{1}-u_{2}\right|<k\right\}}\left(\left|\nabla u_{1}\right|^{p-2} \nabla u_{1}-\left|\nabla u_{2}\right|^{p-2} \nabla u_{2}\right)\left(\nabla u_{1}-\nabla u_{2}\right) d\sigma_{g}
 $$

for $p<2, $ then $ \nabla u_ {1} - \nabla u_{2}= 0 $ a.e. and consequently $ T_{k} \left (u_{2} \right) = T_{k} \left (u_{2} \right) $ for all
$k>0.$ It is clear that $u_{1}-u_{2} =c,$ using the fact that $u_{1}= u_{2}= 0$on$ \partial M, $ then we
concludes that $ u_ {1} = u_ {2} $ a.e. Hence the result. \hfill$\blacksquare$

\subsection{Some generalizations}

The notion of the entropy solution can be defined for a very large class of nonlinear elliptic operators, for example if we consider the next problem
$$-\operatorname{div}(a(x, u, \nabla u))=F(x, u)$$

$\operatorname{with} a(x, s, \xi) : \mathbb{R}^{N} \times \mathbb{R} \times \mathbb{R}^{N} \rightarrow \mathbb{R}^{N}$ is a function of Carathéodory verifying
$(\mathbf{H} \mathbf{1}) :|a(x, s, \xi)| \leq c\left(|\xi|^{p-1}+|s|^{p-1}+k(x)\right), p . p . \quad x \in M,(s, \xi) \in \mathbb{R}^{N} \times \mathbb{R},$

$(\mathbf{H} 2) : a(x, s, \xi) \xi \geq \lambda|\xi|^{p}, p . p . \quad x \in M, \quad(s, \xi) \in \mathbb{R}^{N} \times \mathbb{R},$

$(\mathbf{H} 3) : F$ is a Carateodory function, continuous and decreasing in $ u $ for $ x $ fixed and measurable in $ x $ for $ u $ fixed. Furthermore, $F(x, 0) \in L^{1}(M)$ et $F(x, c) \in L_{l o c}^{1}(M)$ if $c \neq 0$ and $\mathrm{si}$
$$ 
G_{c}(x)=\sup _{|u| \leq c}|F(x, u)|,
 $$

$G_{c} \in L_{l o c}^{1}(M)$ for all $c>0.$\\

So under the conditions $(\mathbf{H} \mathbf{1}),(\mathbf{H} 2)$ et $(\mathbf{H} 3)$ we can define the notion of the solution
in the sense of entropy. Regarding the uniqueness of the solution, in general the result is not
true but if $div(a(x, u, \nabla u ) )=\Delta_{p} u,$ then we can demonstrate the uniqueness of the solution in the sense of entropy.

\vfill
\textbf{Keywords : } Quasi-linear elliptic equations, variational methods, functional spaces, entropy solution, Riemannian manifold, space Marcinkiewicz .


\end{document}